\documentclass[10.5pt]{article}
\usepackage{amsfonts,mathrsfs}

\def\sqr#1#2{{\vcenter{\vbox{\hrule height.#2pt
              \hbox{\vrule width.#2pt height#1pt \kern#1pt \vrule width.#2pt}
              \hrule height.#2pt}}}}
\def\signed #1{{\unskip\nobreak\hfil\penalty50
              \hskip2em\hbox{}\nobreak\hfil#1
              \parfillskip=0pt \finalhyphendemerits=0 \par}}
\def\endpf{\signed {$\sqr69$}}
\def\3n{\negthinspace \negthinspace \negthinspace }
\def\2n{\negthinspace \negthinspace }
\def\1n{\negthinspace }

\def\dbE{\mathbb{E}}
\def\dbF{\mathbb{F}}

\def\dbH{\mathbb{H}}

\def\dbN{\mathbb{N}}

\def\dbP{\mathbb{P}}

\def\dbR{\mathbb{R}}
\def\dbS{\mathbb{S}}

\def\sB{\mathscr{B}}

\def\sE{\mathscr{E}}

\def\sH{\mathscr{H}}

\def\sM{\mathscr{M}}
\def\sN{\mathscr{N}}

\def\sU{\mathscr{U}}

\def\={\buildrel \triangle \over =}

\def\ds{\displaystyle}

\def\ns{\noalign{\ss}}

\def\a{\alpha}
\def\b{\beta}

\def\d{\delta}
\def\e{\varepsilon}
\def\z{\zeta}

\def\si{\sigma}
\def\t{\tau}
\def\f{\varphi}
\def\th{\theta}
\def\o{\omega}

\def\i{\infty}
\def\G{\Gamma}
\def\D{\Delta}
\def\Th{\Theta}
\def\L{\Lambda}

\def\F{\Phi}
\def\Om{\Omega}

\def\cF{{\cal F}}

\def\cH{{\cal H}}

\def\cP{{\cal P}}
\def\cQ{{\cal Q}}

\def\cU{{\cal U}}

\def\cX{{\cal X}}
\def\cY{{\cal Y}}
\def\cZ{{\cal Z}}

\def\BF{{\bf F}}

\def\ss{\smallskip}
\def\ms{\medskip}

\def\q{\quad}
\def\qq{\qquad}
\def\hb{\hbox}

\def\lan{\mathop{\langle}}
\def\ran{\mathop{\rangle}}

\def\esssup{\mathop{\rm esssup}}

\def\h{\widehat}
\def\wt{\widetilde}
\def\cd{\cdot}
\def\cds{\cdots}

\def\ae{\hbox{\rm a.e.{ }}}
\def\as{\hbox{\rm a.s.{ }}}

\def\({\Big (}
\def\){\Big )}
\def\[{\Big[}
\def\]{\Big]}

\def\({\Big (}
\def\){\Big )}
\def\[{\Big[}
\def\]{\Big]}

\def\bde{\begin{definition}}
\def\ede{\end{definition}}
\def\be{\begin{equation}}
\def\bel{\begin{equation}\label}
\def\ee{\end{equation}}
\def\bt{\begin{theorem}}
\def\et{\end{theorem}}
\def\bc{\begin{corollary}}
\def\ec{\end{corollary}}
\def\bl{\begin{lemma}}
\def\el{\end{lemma}}
\def\bp{\begin{proposition}}
\def\ep{\end{proposition}}
\def\bas{\begin{assumption}}
\def\eas{\end{assumption}}
\def\br{\begin{remark}}
\def\er{\end{remark}}
\def\ba{\begin{array}}
\def\ea{\end{array}}
\def\ed{\end{document}}

\def\square#1{\vbox{\hrule\hbox{\vrule height#1%
     \kern#1\vrule}\hrule}}
\def\rectangle#1#2{\vbox{\hrule\hbox{\vrule height#1%
     \kern#2\vrule}\hrule}}

\font\tenbb=msbm10 \font\sevenbb=msbm7 \font\fivebb=msbm5

\newfam\bbfam
\scriptscriptfont\bbfam=\fivebb \textfont\bbfam=\tenbb
\scriptfont\bbfam=\sevenbb

\newtheorem{lemma}{Lemma}[section]
\newtheorem{remark}{Remark}[section]

\newtheorem{theorem}{Theorem}[section]
\newtheorem{corollary}{Corollary}[section]

\newtheorem{definition}{Definition}[section]
\newtheorem{proposition}{Proposition}[section]
\newtheorem{assumption}{Assumption}[section]

\makeatletter
   
   \@addtoreset{equation}{section}
\makeatother

\usepackage{amsfonts}
\usepackage{latexsym}
\usepackage{amsmath}
\usepackage{amssymb}
\usepackage{color}

\def\ges{\geqslant}
\def\les{\leqslant}
\def\rf{\eqref}

\begin{document}

\title{Spike Variations for\\ Stochastic Volterra Integral Equations}

\author{ Tianxiao Wang \footnote{School of
Mathematics, Sichuan University, Chengdu, P. R. China. Email: wtxiao2014@scu.edu.cn.
This author was supported by National Natural Science Foundation of China (No. 11971332
and 11931011) and the Science Development Project of Sichuan University under grant
2020SCUNL201.
} ~~ and ~~~Jiongmin Yong \footnote{Department of
Mathematics, University of Central Florida, Orlando, USA. Email: Jiongmin.Yong@ucf.edu. This author is supported in part by NSF grant DMS-1812921.}
}

\maketitle

\begin{abstract}

Spike variation technique plays a crucial role in deriving Pontryagin's type maximum principle of optimal controls for differential equations of several types, including ordinary differential equations (ODEs), partial differential equations (PDEs), and stochastic differentia equations (SDEs), when the control domains are not assumed to be convex. This technique also applies to (deterministic forward) Volterra intrgral equations (FVIEs). It is natural to expect that such a technique could be extended to the case of (forward) stochastic Volterra integral equations (FSVIEs). However, by mimicking the case of SDEs, one encounters an essential difficulty of handling an involved quadratic term. To overcome the difficulty, we introduce an auxiliary process for which one can use It\^o's formula, and adopt a trick used in linear-quadratic stochastic optimal control problems. Then a suitable representation of the above-mentioned quadratic form is obtained, and the second order adjoint equations are derived. Consequently, the maximum principle of Pontryagin type is established. Some relevant extensions are investigated as well.
\end{abstract}

\ms

\bf Keywords. \rm  Backward stochastic Volterra integral equations, spike variation, second-order adjoint equations, non-convex control domain.

\ms

\bf AMS Mathematics subject classification. \rm 45D05, 60H20, 93E20

\section{Introduction}

We begin with a main motivation of our study. Let $T>0$ and $n,m\in\dbN$ be given,
$(\Om,\cF,\dbP)$ be a complete probability space, on which a one-dimensional
standard Brownian process $W(\cd)$ is defined, with the completed augmented natural filtration denoted by $\dbF=\{\cF_t\}_{0\les t\les T}$. Consider the following controlled (forward) stochastic Volterra integral equation (FSVIE, for short):
\bel{state}X(t)=\f(t)+\int_0^tb(t,s,X(s),u(s))ds+\int_0^t\si(t,s,X(s),u(s))dW(s),\q t\in[0,T],\ee
with the cost functional
\bel{cost3.1}J(u(\cd))=\dbE\Big[h(X(T))+\int_0^Tg(s,X(s),u(s))ds\Big].\ee
Here $u(\cd)$ is a {\it control process} taking values in a non-empty measurable set $U\subseteq\dbR^m$, and $X(\cd)$ is the corresponding {\it state process} valued in $\dbR^n$, $\f:[0,T]\times\Om\to\dbR^n$ is called a {\it free term}, $(b,\si):[0,T]^2\times\Om\times\dbR^n\times\dbR^m\to
\dbR^n\times\dbR^n$ is called the {\it generator} of the controlled FSVIE, $h:\Om\times\dbR^n\to\dbR$ $g:[0,T]\times\Om\times\dbR^n\times\dbR^m\to\dbR$ are called the {\it terminal cost} and {\it running cost rate}, respectively. Under certain conditions, for any $u(\cd)\in\sU_{ad}$ (a set of all admissible controls, to be defined later), \rf{state} admits a unique solution $X(\cd)$ such that the cost functionl $J(u(\cd))$ is well-defined.
The optimal control problem can be stated as follows:

\ss

\bf Problem (C). \rm Find a $\bar u(\cd)\in\sU_{ad}$ such that
\bel{J=inf}J(\bar u(\cd))=\inf_{u(\cd)\in\sU_{ad}}J(u(\cd)).\ee

If the above is true, $\bar u(\cd)$ is called an {\it open-loop optimal control}, the corresponding state process $\bar X(\cd)$ and $(\bar X(\cd),\bar u(\cd))$ are called the ({\it open-loop}) {\it optimal state process} and an ({\it open-loop}) {\it optimal pair}, respectively.

\ss

Pontryagin's maximum principle for (deterministic) controlled ordinary differential equations (ODEs) was formulated and proved by Boltyanskii--Gamkrelidze--Pontryagin in 1956 \cite{BGP-1956} (see also \cite{PBGM-1962}). Since then, the research in this direction
attracted many authors. Shortly after, people attempted to extend the
results to stochastic differential equations (SDEs, for short).  In 1965, Kushner firstly proved a maximum principle for SDEs with deterministic coefficients, for which the diffusion is independent of the state and control (\cite{Kushner-1965}). Therefore, the diffusion of the variational system vanishes so that the adjoint equation is an ODE, whose solution is the adjoint state process. In 1972, Kushner further considered the case that the diffusion of the state equation contains the state (\cite{Kushner-1972}), in which the adjoint process was not successfully characterized. In 1973, Bismut \cite{Bismut-1978} introduced the duality between an initial value problem of linear SDE (called FSDE) and a terminal value problem of a linear SDE (now called BSDE, a name coined by Pardoux--Peng in 1990  \cite{Pardoux-Peng-1990}). With this, he was able to prove the maximum principle for FSDEs with the diffusion contains state and control, provided the control domain is convex; the adjoint process was successfully characterized by the adapted solution (again, named by Pardoux--Peng in 1990) of corresponding linear BSDE. There were a number of follow-up works. Among them, we mention Haussmann (\cite{Haussmann-1976}) and Bensoussan (\cite{Bensoussan-1981}). The maximum principle for general controlled SDEs  of the following form
\bel{FSDE}X(t)=x+\int_0^tb(s,X(s),u(s))ds+\int_0^t\si(s,X(s),u(s))dW(s),\q t\in[0,T],\ee%
without assuming the convexity of the control domain was proved by Peng in 1990 (\cite{Peng-1990}), with the development of spike variation technique for FSDEs. We refer to the monograph of Yong--Zhou \cite{Yong-Zhou-1999} for a self-contained presentation.
By the way, spike variation techniques for partial differential equations (PDEs, for short) and abstract evolution equations in Hilbert/Banach spaces were also fully developed during 1970--1990s. See Li--Yong \cite{Li-Yong-1995}, Casas--Yong \cite{Casas-Yong-1995}, Hu--Yong \cite{Hu-Yong-1995}, Fattorini \cite{Fattorini-1999}, and so on, for some systematic presentations.

\ss

Integral equation was firstly introduced by N. H. Abel in 1825. The following-up contributors include Liouville, Fredholm, Hilbert, Wiener, Bellman, etc. We refer to B\^ocher \cite{Bocher-1909}, Yosida \cite{Yosida-1960} for two early monographs in this topic. (Forward) Volterra integral equation (FVIEs, for short) was introduced by Italian mathematician Volterra. In contrast with ODEs, FVIEs bring us new theories (e.g., Volterra operator theory in Gohberg--Krein \cite{Gohberg-Krein-book-2004}) and new phenomena (e.g., the hereditary property in Bellman--Cooke \cite{Bellman-Cooke-1963}, Volterra \cite{Volterra-book-1930}). In the real world, there are many models that cannot be simply described by ODEs but by (deterministic) FVIEs. Therefore, up until now, it is still an active research topic, see, e.g., the recent monograph of Brunner \cite{Brunner-2017}. For the optimal control theory of FVIEs, to our best knowledge, the earliest work was due to Friedman \cite{Friedman-1964} in 1964. Later in 1967, Vinokurov \cite{Vinokurov-1969} investigated the optimal control for FVIEs with various type constraints and closed control domains (described by level sets of smooth functions). There are quite a few literature concerning optimal control of FVIEs in the past several decays. We mention some of them here: Kamien--Muller \cite{Kamien-Muller-1976}, Medhin \cite{Medhin-1986}, Belbas \cite{Belbas-2008}, and Bonnans--de la Vega--Dupuis \cite{Bonnans-de la Vega-Dupuis-2013}. Recently, Lin--Yong \cite{Lin-Yong-2018} established a maximum principle for singular FVIEs with some state constraint by spike variation technique.

\ss

In contrast with the existing literatures for optimal control of (deterministic) FVIEs, there were much less literature on optimal control of FSVIEs, to our best knowledge. The simple reason is that the theory of adjoint equation for FSVIEs, i.e., the backward SVIEs (BSVIEs, for short) was not available until the appearance of the work on Type-I BSVIEs by Lin (\cite{Lin-2002}) in 2002. Even then, the result was not connected with the maximum principle of optimal control for FSVIEs. In 2006, Yong introduced Type-II BSVIEs and established maximum principle for FSVIEs the first time (\cite{Yong-2006}), provided the control domain is convex (see also Yong \cite{Yong-2008}). There are some further extensions of maximum principle for FSVIEs appeared lately. Here are some of those results. Shi--Wang--Yong \cite{Shi-Wang-Yong-2015} studied the case of forward-backward SVIEs. Agram--{\O}ksendal \cite{Agram-Oksendal-2015} presented some investigations by means of Malliavin calculus, when the control domain $U$ is open. Wang--Zhang \cite{Wang-Zhang-2017} considered the case with $U$ being a closed set with the necessary optimality conditions were obtained from those directions that approximate convex perturbations at the optimal control are admissible. In such a case, spike variations are not necessary.

\ss

Recently, Wang \cite{Wang-2020} studied the general case without assuming the convexity of the control domain $U$ and attempted to develop the spike variation techniques for FSVIEs. In doing so, the following form {\it quadratic functional}
\bel{sE(e)0}\sE(\e)=\dbE\Big[X_1^\e(T)^\top h_{xx}X_1^\e(T)+\int_0^TX_1^\e(t)^\top H_{xx}(t)X_1^\e(t) dt\Big],\ee
appeared in the variation of the cost functional, where $h_{xx}$ and $H_{xx}(\cd)$ are $\dbS^n$ (the set of all real $(n\times n)$ symmetric matrices) valued functions only depend on the optimal pair $(\bar X(\cd),\bar u(\cd))$, and $X_1^\e(\cd)$ is the solution of the first-order variational equation, depending on the spike variation of the control. In the case that the control domain is convex, convex variation of the control is allowed and by doing that the above quadratic form does not appear. Thus, it is not necessary for \cite{Agram-Oksendal-2015, Shi-Wang-Yong-2015, Wang-Zhang-2017, Yong-2006, Yong-2008} to introduce techniques of handling the above. On the other hand, recall that similar quadratic form as $\sE(\e)$ appears in the FSDE case, which can be handled by a duality between the second-order adjoint SDE and the variational SDE satisfied by $X_1^\e(\cd)X_1^\e(\cd)^\top$ (see Peng \cite{Peng-1990}). It is unlikely that for FSVIEs, one can obtain a good-looking linear FSVIE for $X_1^\e(\cd)X_1^\e(\cd)^\top$. Hence, the idea of Peng \cite{Peng-1990} cannot be directly borrowed here. In fact, in Wang \cite{Wang-2020}, the author had to introduce two abstract and implicit operator-valued stochastic processes, which are regarded as a replacement of the second-order adjoint processes, in the statement of the maximum principle. Moreover, there was no second-order adjoint equation derived in \cite{Wang-2020}.

\ss

In the current paper, we are going to develop a spike variation technique for FSVIEs. The main ideas can be described as follows. Recall that in the study of linear-quadratic (LQ, for short) optimal controls, people could derive the Riccati equation by applying the It\^o's formula to the map $t\mapsto X(t)^\top P(t)X(t)$ with $P(\cd)$ being a symmetric matrix valued process which characterizes the quadratic value function of the problem. It turns out that this idea can be adopted to our current problem. Together with the introduction of a suitable auxiliary process, we are able to represent the quadratic form involved $X_1^\e(\cd)$ in a desired fashion. Consequently, the second-order adjoint equation will be derived leading to the maximum principle. Moreover, we will go a little further. Recall that in 2007, Mou--Yong (\cite{Mou-Yong-2007}) established a variational formula for an SDE state equation and the possibly vector-valued cost/payoff functionals. Those results are applicable to multi-objective problems and multi-person differential games of SDEs. The proof of such a variational formula (for SDEs) is almost the same as that for single-objective optimal control problem. Inspired by this, in the current paper, for FSVIEs, we will do the same thing, so that it will lead to necessary conditions for the Nash equilibria of multi-person dynamic games for FSVIEs. By letting the number of players go to infinite, and under some structure conditions, it is expected to obtain the necessary conditions for the equilibria of the mean-field games governed by FSVIEs, which includes the standard SDE case. We hope to report more details in our forthcoming publications.

\ss

The rest of this paper is organized as follows. In Section 2, we introduce some preliminary notations and assumptions. In Section 3, we present the variations of the state and the cost functional under the spike variation of the control process, whose detailed and lengthy proof will be given in Section 7. In Section 4, we give a representation of the involved  quadratic functional by introducing a system of BSVIEs, called the {\it second order adjoint equations}. In Section 5, we prove the well-posedness of the aforementioned second order adjoint equations. We present the Pontryagin's type maximum principle and some extensions to multi-person dynamic games for FSVIEs in Section 6. Finally, some concluding remarks are collected in Section 8.

\section{Preliminaries}

Let us introduce some basic spaces. First of all, let $\dbR^n$ and $\dbR^{n\times d}$ be the usual $n$-dimensional space of real numbers and the set of all $(n\times d)$ real matrices, respectively. Also, let $\dbS^n$ be the set of all $(n\times n)$ real symmetric matrices. Next, for any Euclidean space $\dbH$ which could be $\dbR^n,\dbR^{n\times d}$, $\dbS^n$, etc., we define:
$$L^p_{\cF_t}(\Omega;\dbH)=\Big\{\xi:\Om\to\dbH\bigm|\xi\hb{ is $\cF_t$-measurable, }\|\xi\|_p\equiv\big(\dbE|\xi|^p\big)^{1\over p}<\infty\Big\},\qq p\in[1,\infty).$$
In an obvious way, we can define $L^\infty_{\cF_t}(\Om;\dbH)$. It is known that for each $p\in[1,\i]$, $L^p_{\cF_t}(\Om;\dbH)$ is a Banach space under the norm $\|\cd\|_p$. When the range space $\dbH$ is clear from the context and is not necessarily to be emphasized, we will omit $\dbH$. In particular, we will denote $L^p_{\cF_T}(\Omega)=L^p_{\cF_T}(\Omega;\dbH)$.

\ss

Next, we introduce spaces of stochastic processes. In order to avoid repetition, all processes $(t,\o)\mapsto\f(t,\o)$ are assumed to be at least $\sB[0,T]\otimes\cF_T$-measurable without further mentioning, where $\sB[0,T]$ is the Borel $\si$-field of $[0,T]$. For $p,q\in[1,\infty)$, $\t\in[0,T)$,
$$\ba{ll}
\ds L^p(\Om;L^q(\t,T;\dbH))\1n=\1n\Big\{\f:[\t,T]\1n\times\1n\Om\to\dbH\bigm|
\dbE\(\int_\t^T\3n|\f(t)|^qdt\)^{p\over q}\3n<\infty\Big\},\\
\ns\ds L^p(\Om;L^\infty(\t,\1n T;\dbH))\1n=\1n\Big\{\f:[\t,T]\1n\times\1n\Om\to\dbH\bigm|
\dbE\(\esssup_{t\in[\t,T]}|\f(t)|^p\)\1n<\1n\infty\Big\},\\
\ns\ds L^p(\Om;C([\t,\1n T];\dbH))\1n=\1n\Big\{\f:[\t,T]\1n\times\1n\Om\to\dbH\bigm|
t\mapsto\f(t,\o)\hb{ is continuous, }\dbE\(\sup_{t\in[\t,T]}|\f(t)|^p\)\1n<\1n\infty\Big\},\\
\ds L^q(\t,T;L^p(\Om;\dbH))\1n=\1n\Big\{\f:[\t,T]\1n\times\1n\Om\to\dbH
\bigm|
\int_\t^T\(\dbE|\f(t)|^p\)^{q\over p}dt<\infty\Big\},\\
\ns\ds L^\infty(\t,T;L^p(\Om;\dbH))\1n=\1n\Big\{\f:[\t,T]\1n\times\1n\Om\to
\dbH\bigm|
\esssup_{t\in[\t,T]}\(\dbE|\f(t)|^p\)^{1 \over p}
<\1n\infty\Big\},\\
\ns\ds C([\t,T];L^p(\Om;\dbH))\1n=\1n\Big\{\f:[\t,T]\1n\times\1n\Om\to\dbH
\bigm|t\mapsto\f(t,\cd)\hb{ is continuous, }\sup_{t\in[0,T]}\(\dbE|\f(t)|^p\)^{1\over p}<\1n\infty\Big\}.\ea$$
The spaces $L^\infty(\Om;L^\infty(\t,T;\dbH))$, $L^\infty(\t,T;L^\infty(\Om;\dbH))$, $L^\infty(\Om;C([\t,T];\dbH))$, and $C([\t,T];$ $L^\infty(\Om;$ $\dbH))$ can be defined in an obvious way. For all $p\in[1,\infty)$, we denote
$$L^p(\t,T;\dbH)=L^p(\t,T;L^p(\Om;\dbH))=L^p(\Om;L^p(\t,T;\dbH)).$$
Note that processes in the above spaces are not necessarily $\dbF$-adapted. The subset of, say, $L^p(\Omega;L^q(\t,T;\dbH))$, consisting of all $\dbF$-progressively measurable processes, is denoted by $L^p_\dbF(\Omega;$ $L^q(\t,T;\dbH))$. All the $\dbF$-progressive version of the other spaces listed in the above can be denoted similarly, by putting $\dbF$ as a subscript; for example, $C_\dbF([\t,T];L^p(\Omega;\dbH))$, and so on. Finally, in the above, when the range space $\dbH$ is clear from the context, we will omit $\dbH$; for example, $L^p_\dbF(\Omega;C([\t,T]))$.

\ss

Next, the upper and lower triangle domains are defined by the following:
$$\ba{ll}
\ns\ds\D^*[\t,T]=\big\{(r,s)\in[\t,T]^2\bigm|\t\les r<s\les T\big\},\\
\ns\ds\D_*[\t,T]=\big\{(r,s)\in[\t,T]^2\bigm|\t\les s<r\les T\big\},\ea\qq\t\in[0,T).$$
We introduce the following spaces.
$$\ba{ll}
\ns\ds L^2_\dbF(\Om;L^2(\D^*[\t,T]))=\Big\{\z:\D^*[\t,T]\times\Omega\to\dbH\bigm|\ae r\in[\t,T],\, \z(r,\cd)\in L^2_\dbF(\Omega;L^2(r,T;\dbH)),\\
\ns\ds\qq\qq\qq\qq\qq\qq\qq\qq\qq\qq\qq\dbE\int_\t^T\int_r^T|\z(r,s)|^2ds dr<\infty\Big\},\\
\ns\ds L^2_\dbF(\Om;L^2(\D_*[\t,T]))=\Big\{\z:\D_*[\t,T]\times\Omega\to\dbH\bigm|\ae r\in[\t,T],\, \z(r,\cd)\in L^2_\dbF(\Omega;L^2(\t,r;\dbH)),\\
\ns\ds\qq\qq\qq\qq\qq\qq\qq\qq\qq\qq\qq\dbE\int_\t^T\int_\t^r|\z(r,s)|^2ds dr<\infty\Big\}.
%
%
%
\ea$$
%

In what follows, for simplicity of the presentation, we will let $K$ be a generic constant which could be different from line to line.
For FSVIE (\ref{state}), we introduce the following assumptions.

\medskip

{\bf(H1)} Let $T>0$ and $\varnothing\ne U\subset\dbR^m$ be measurable. Let $b,\si:\D_*[0,T]\times\dbR^n\times U\times\Omega\to\dbR^n$ be measurable such that
$s\mapsto\big(b(t,s,x,u),\si(t,s,x,u)\big)$
is $\dbF$-progressively measurable on $[0,t]$, $x\mapsto(b(t,s,x,u),\si(t,s,x,u))$ is twice continuously differentiable with uniformly bounded first and second order derivatives, and for some constant $L>0$,
$$|b(t,s,0,u)|+|\si(t,s,0,u)|\les L,\qq(t,s,u)\in\D_*[0,T]\times U.$$
Further, for $f=b,\si,b_x,\si_x$, the map $(t,u)\mapsto f(t,s,x,u)$ is continuous
uniformly in all the arguments, and for $f=b,\si$, the map $(x,u)\mapsto f_{xx}(t,s,x,u)$ is continuous uniformly in all arguments.

\ms

Next, we introduce the set of all admissible controls
$$\sU_{ad}\=\Big\{u:[0,T]\times\Om\to U\bigm|u(\cd)\hb{ is $\dbF$-progressively measurable}\Big\}.$$
Now, we present a result of well-posedness for our state equation, whose proof is standard.

\bp\label{FSVIE} \sl Let {\rm(H1)} hold, $\f(\cd)\in C_\dbF([0,T];L^p(\Om;\dbR^n))$, for some $p\ges1$. Then for any $u(\cd)\in\sU_{ad}$, \rf{state} admits a unique solution $X(\cd)\equiv X(\cd\,;u(\cd))\in C_{\dbF}([0,T];L^p(\Om;\dbR^n))$ such that
\bel{|X|}\sup_{t\in[0,T]}\1n\dbE|X(t)|^p\1n\les\1n K\2n\sup_{t\in[0,T]}\1n\dbE\Big[|\f(t)|^p\1n+\1n\(\1n\int_0^t\3n|b(t,s,0,u(s))|ds\)^p
\3n+\1n\(\1n\int_0^t\3n|\si(t,s,0,u(s))|^2ds\)^{p\over2}\Big].\ee
Hereafter, $K$ will be a generic constant which could be different from line to line. Moreover, if $(\bar X(\cd),\bar u(\cd))$ is another state-control
pair. Then the following stability estimate holds:
\bel{|X-X|*}\ba{ll}
\ns\ds\sup_{t\in[0,T]}\1n\dbE|X(t)-\bar X(t)|^p\1n\les\1n K\2n\sup_{t\in[0,T]}\1n\dbE\Big[\(\1n\int_0^t\3n|b(t,s,\bar X(s),u(s))-b(t,s,\bar X(s),\bar u(s))|ds\)^p\\
\ns\ds\qq\qq\qq\qq\qq\qq+\(\1n\int_0^t\3n|\si(t,s,\bar X(s),u(s))-\si(t,s,\bar X(s),\bar u(s))|^2ds\)^{p\over2}\Big].\ea\ee

\ep

In what follows, we will fix $p\ges4$. For the involved functions $h$ and $g$ in (\ref{cost3.1}), we make the following assumption.

\ms

{\bf (H2)} Let $h:\dbR^n\times\Om\to\dbR$, $g:[0,T]\times\dbR^n\times U\times\Om\to\dbR$
be measurable such that $x\mapsto(h(x),g(t,x,u))$ is twice continuously differentiable with
$$|g_x(t,0,u)|\les L,\q|h_{xx}(x)|+|g_{xx}(t,x,u)|\les L,\qq(t,x,u)\in[0,T]\times\dbR^n\times U,$$
for some constant $L>0$ and $(t,x,u)\mapsto g_{xx}(t,x,u)$ is uniformly continuous.

\ms

We see that under (H1)--(H2), Problem (C) is well-formulated. Further, if we let $h$ and $g$ take values in $\dbR^\ell$, then our study covers multi-objective problems as well as multi-person non-zero sum dynamic games governed by FSVIEs (if $\ell\ges2$). Finally, we point out that (H1)--(H2) can be relaxed somehow, for example, we may allow $b,\si$, as well as $g$ to have some growth in $u$. We prefer not to get into the most generality to keep the presentation simple.

\section{Variations of the State and the Cost Functional }

In this section, we are going to present the variations of the state process and the cost functional under the spike variation of the control process. For Problem (C), such a result will lead to the first order optimality condition in the form of variational inequality.

\ss

To state the main result of this section, let $(\bar X(\cd),\bar u(\cd))$ be any fixed state-control pair of the state equation, which could be optimal, in particular. Choose any $u\in U$, $\t\in[0,T)$, and let $\e>0$ be sufficiently small such that $\t+\e\les T$, we define
\bel{u^e}u^{\e}(\cd)= u{\bf1}_{[\t,\t+\e]}(\cd)+\bar u(\cd){\bf1}_{[0,T]\setminus[\t,\t+\e]}(\cd).\ee
This is called a {\it spike variation} of $\bar u(\cd)$ in the direction of $u$ on $[\t,\t+\e]$. Let $X^\e(\cd)$ be the state process corresponding to $u^\e(\cd)$. Our goal is to obtain expansions of $X^\e(\cd)-\bar X(\cd)$ and $J(u^\e(\cd))-J(\bar u(\cd))$ in terms of the powers of $\e$.

\ss

To begin, we introduce the following abbreviations:
$$\ba{ll}
\ns\ds b_x(t,s)=b_x(t,s,\bar X(s),\bar u(s)),\qq b^i_{xx}(t,s)=b^i_{xx}(t,s,\bar X(s),\bar u(s)),\qq1\les i\les n,\\
\ns\ds\si_x(t,s)=\si_x(t,s,\bar X(s),\bar u(s)),\qq\si^i_{xx}(t,s)=\si^i_{xx}(t,s,\bar X(s),\bar u(s)),\qq1\les i\les n,\\
\ns\ds g_x(s)=g_x(s,\bar X(s),\bar u(s)),\qq g_{xx}(s)=g_{xx}(s,\bar X(s),\bar u(s)),\\
\ns\ds h_x=h_x(\bar X(T)),\qq h_{xx}=h_{xx}(\bar X(T)),\\
\ns\ds\d b(t,s)=b(t,s,\bar X(s),u)-b(t,s,\bar X(s),\bar u(s)),\\
\ns\ds\d\si(t,s)=\si(t,s,\bar X(s),u)-\si(t,s,\bar X(s),\bar u(s)),\\
\ns\ds\d\si_x(t,s)=\si_x(t,s,\bar X(s),u)-\si_x(t,s,\bar X(s),\bar u(s)).\ea$$
With the above notations, we introduce the following variational equations:
\bel{X_1-equation} X_1^\e(t)\1n=\2n\int_0^t\3n b_x(t,s)X_1^\e(s)ds+\2n\int_0^t\2n\big[\si_x(t,s)X_1^\e(s)+\d\si(t,s){\bf1}_{[\t,\t+\e]}(s)
\big]dW(s),\q t\in[0,T],\ee
and
\bel{X_2-equation}\ba{ll}
\ss\ds X_2^\e(t)=\int_0^t\(b_x(t,s)X_2^\e(s)+{1\over2}b_{xx}(t,s)X_1^\e(s)^2+\d b(t,s){\bf1}_{[\t,\t+\e]}(s)\)ds\\
\ss\ds\qq\qq+\int_0^t\(\si_x(t,s)X_2^\e(s)+{1\over2}
\si_{xx}(t,s)X_1^\e(s)^2+\d\si_x(t,s)X_1^\e(s){\bf1}_{[\t,\t+\e]}(s)\)dW(s),\\
\ns\ds\qq\qq\qq\qq\qq\qq\qq\qq\qq\qq\qq\qq\qq\qq\qq t\in[0,T],\ea\ee
where
$$b_{xx}(t,s)X_1^\e(s)^2=\begin{pmatrix}X_1^\e(s)^\top b^1_{xx}(t,s)X_1^\e(s)\\ \vdots\\ X_1^\e(s)^\top b^n_{xx}(t,s)X_1^\e(s)\end{pmatrix},$$
and $\si_{xx}(t,s)X_1^\e(s)^2$ is similar. Under (H1), the above variational equations \rf{X_1-equation} and \rf{X_2-equation} have unique strong solutions $X_1^\e(\cd)$ and $X_2^\e(\cd)$. We also introduce the following first order adjoint equations:
\bel{1st-adjoint-equation}\left\{\2n\ba{ll}
\ns\ds\eta(t)=h_x^\top-\int_t^T\z(s)dW(s),\qq t\in[0,T],\\
\ns\ds Y(t)=g_x(t)^\top+b_x(T,t)^\top h_x^\top+\si_x(T,t)^{\top}\z(t)\\
\ns\ds\qq\q+\int_t^T\(b_x(s,t)^\top Y(s)+\si_x(s,t)^\top Z(s,t)\)ds-\int_t^TZ(t,s)dW(s),\q t\in[0,T].\ea\right.\ee
The first equation of the above is the simplest BSDE which is just a martingale representation. The second is a Type-II BSVIE which has a unique adapted M-solution (\cite{Yong-2008}). It is clear that the solutions $(\eta(\cd),\z(\cd))$ and $(Y(\cd),Z(\cd\,,\cd))$ of \rf{1st-adjoint-equation} are only depending on the pair $(\bar X(\cd),\bar u(\cd))$ (and independent of the spike variation of the control). Finally, let us define the following Hamiltonian:
\bel{H}\ba{ll}
\ns\ds H(s,x,u,\eta(s),\z(s),Y(\cd),Z(\cd\,,s))\1n=\1n\lan\eta(s),b(T,s,x,u)\ran
\1n+\1n\lan\z(s),\si(T,s,x,u)\ran\1n+g(s,x,u)\\
\ns\ds\qq\qq\qq\qq\qq\qq\qq\;+\dbE_s\2n\int_s^T\3n\(\2n\lan Y(t),b(t,s,x,u)\ran\1n+\1n\lan Z(t,s),\si(t,s,x,u)\ran\2n\)dt.\ea\ee
Hereafter, $\dbE_s=\dbE[\,\cd\,|\cF_s]$ is the conditional expectation operator. Then the main result of this section can be stated as follows.

\bt\label{variation} \sl Let {\rm(H1)--(H2)} hold. Then for any $k\ges1$,
\bel{|X_1|,|X_2|}\left\{\2n\ba{ll}
\ns\ds\sup_{s\in[0,T]}\dbE|X^\e(s)-\bar X(s)|^k=O(\e^{k\over2}),\qq\sup_{s\in[0,T]}\dbE|X_1^\e(s)|^k=O(\e^{k\over2}),\\
\ns\ds\sup_{s\in[0,T]}\dbE|X^\e(s)-\bar X(s)-X_1^\e(s)|^k=O(\e^k),\qq\sup_{s\in[0,T]}\dbE|X_2^\e(s)|^k=O(\e^k),\\
\ns\ds\sup_{s\in[0,T]}\dbE|X^\e(s)-\bar X(s)-X_1^\e(s)-X_2^\e(s)|^2=o(\e^2).\ea\right.\ee
Moreover,
\bel{J-J}J(u^\e(\cd))-J(\bar u(\cd))=\dbE\int_\t^{\t+\e}\D H(s)ds+{1\over2}\sE(\e)+o(\e),\ee
where
\bel{DH(e)}\ba{ll}
\ns\ds \D H(s)=H(s,\bar X(s),u,\eta(s),\z(s),Y(\cd),Z(\cd\,,s))\\
\ns\ds\qq\qq\qq\qq\qq\q-H(s,\bar X(s),\bar u(s),\eta(s),\z(s),Y(\cd),Z(\cd\,,s)),\ea\ee
and
\bel{sE(e)}\ba{ll}
\ns\ds\sE(\e)=\dbE\Big[X_1^\e(T)^\top h_{xx}(\bar X(T))X_1^\e(T)\\
\ns\ds\qq\qq\qq\q+\int_0^T\3n X_1^\e(s)^\top H_{xx}(s,\bar X(s),\bar u(s),\eta(s),\z(s),Y(\cd),Z(\cd\,,s))X_1^\e(s)ds\Big]\\
\ns\ds\qq\equiv\dbE\Big[X_1^\e(T)^\top h_{xx}X_1^\e(T)+\int_0^T\3n
X_1^\e(s)^\top H_{xx}(s)X_1^\e(s)ds\Big].\ea\ee
\et

The above is essentially shown in \cite{Wang-2020} with the main idea of \cite{Peng-1990} (see also \cite{Yong-Zhou-1999}). For convenience of the readers, we provide a detailed proof in Section 7. We point out that in the last estimate of \rf{|X_1|,|X_2|}, $2$ cannot be replaced by ``any $k\ges1$''.

\section{Representation of the Quadratic Form.}

Recalling \rf{X_1-equation} (suppressing the superscript $\e$ in $X_1^\e(\cd)$), we have
$$\ba{ll}
\ns\ds X_1(t)=0,\qq\qq t\in[0,\t],\\
\ns\ds X_1(t)= \int_\t^tb_x(t,s)X_1(s)ds+\int_\t^t \big[\si_x(t,s)X_1(s)+\d\si^\e(t,s)\big]dW(s),\q t\in[\t,T],\ea$$
where
$$\d\si^\e(t,s)\equiv\d\si(t,s){\bf1}_{[\t,\t+\e]}(s).$$
Thus, $X_1(\cd)$ is the image of $\d\si^\e(\cd\,,\cd)$ under a linear operator, i.e.,
$$X_1(t)=\F[\d\si^\e(\cd\,,\cd)](t),\qq\qq t\in[\t,T],$$
for some operator $\F$. Let
$$X_1(T)=\F[\d\si^\e(\cd\,,\cd)](T)\equiv\h\F[\d\si^\e(\cd\,,\cd)].$$
Consequently, by misuse the notation, we roughly have
\bel{sE(e)**}\sE(\e)=\lan\big(\h\F^*h_{xx}\h\F+\F^*H_{xx}\F\big)\d\si^\e,\d\si^\e\ran.\ee
In the above, the dependence of $\sE(\e)$ on $\d\si^\e(\cd\,,\cd)$ is not directly enough, due to the fact that the operator $\F$ is complicated. Therefore, \rf{sE(e)**} is not an integral of a bilinear form in $\dbR^n$.
The purpose of this section is to obtain a more direct form of $\sE(\cd)$ in terms of $\d\si^\e(\cd\,,\cd)$, whose procedure leads to the second adjoint equation.
To express the main idea, let us first look at the following simple calculation. Such an idea comes from the standard study of LQ problems.

\ms

Consider the following linear SDE:
$$\left\{\2n\ba{ll}
\ns\ds dX(s)=\big[AX(s)+b(s)\big]ds+\big[CX(s)+\si(s)\big]dW(s),\qq s\in[t,T],\\
\ns\ds X(t)=x,\ea\right.$$
with suitable matrices $A,C$ and processes $b(\cd),\si(\cd)$. Also, we are given the following quadratic functional
$$J=\dbE\Big[X(T)^\top GX(T)+\int_t^T\(X(s)^\top QX(s)+2q(s)^\top X(s)\)ds\Big],$$
with suitable $G,Q$, and $q(\cd)$. We want to absorb the first term on the right-hand side of the above into the second one. To this end, we let $(P(\cd),\L(\cd))$ be the adapted solution to the following BSDE (for symmetric matrix valued processes):
$$\left\{\2n\ba{ll}
\ns\ds dP(t)=-\G(t)dt+\L(t)dW(t),\qq t\in[0,T],\\
\ns\ds P(T)=G,\ea\right.$$
with $\G(\cd)$ undetermined (symmetric matrix valued process). Then by It\^o's formula (suppressing $s$), we have
$$d(PX)=\(-\G X+P(AX+b)+\L(CX+\si)\)ds+\(\L X+P(CX+\si)\)dW.$$
Therefore,
$$\ba{ll}
\ns\ds d(X^\top PX)=\Big[(X^\top A^\top+b^\top)PX+X^\top(-\G X+PAX+Pb+\L CX+\L\si)\\
\ns\ns\qq\qq\qq\qq+(X^\top C^\top+\si^\top)(\L X+PCX+P\si)\Big]ds+\{\cds\}dW(s)\\
\ns\ds=\Big[X^\top\(-\G+PA+A^\top P+C^\top PC+\L C+C^\top\L\)X^\top+X^\top\(Pb+C^\top P\si+\L\si)\\
\ns\ds\qq+\(b^\top P+\si^\top PC+D^\top\L\)X+\si^\top P\si\Big]ds+\{\cds\}dW(s). \ea$$
Then
$$\ba{ll}
\ns\ds J=\dbE\Big[X(T)^\top GX(T)+\int_t^T\(X(s)^\top QX(s)+2q(s)^\top X(s)\)ds\Big]\\
\ns\ds\q=\dbE\int_t^T\Big[X(s)^\top\(-\G+PA+A^\top P+C^\top PC+\L C+C^\top\L+Q\)X(s)\\
\ns\ds\qq\qq+2\(b^\top P+\si^\top PC+\si^\top\L+q^\top\)X+\si^\top P\si\Big]ds+x^{\top}P(t)x.\ea$$
We see that the term $X(T)^\top GX(T)$ has been absorbed. Further, it is easy to see that to elimnate the term $X(s)^\top(\cds)X(s)$ under the integral, we may choose
$$\G=PA+A^\top P+C^\top PC+\L C+C^\top\L+Q,$$
so that one has
$$J=x^{\top}P(t)x+\dbE\int_t^T\Big[2\(Pb+C^\top P\si+\L\si+q\)^\top X+\si^\top P\si\Big]ds.$$
Such an idea will play a crucial role below.

\ms

Note that the above could not apply directly to $\sE(\e)$ as $t\mapsto X_1(t)$ does not satisfy an SDE, and therefore, the It\^o's formula could not apply. To overcome this difficulty, we introduce an auxiliary process $\cX_1(\cd\,,\cd)$ as follows
\bel{Additional-X-1}\ba{ll}
\ns\ds\cX_1(t,r)\1n=\2n\int_\t^r\3n b_x(t,s)X_1(s)ds\1n+\2n\int_\t^r\3n\big[\si_x(t,s)X_1(s)
+\d\si^\e(t,s)\big]dW(s),\q\t\les r\les t\les T.\ea\ee
Similar form also appeared in \cite{Wang-2016-submitted}. Note that $r\mapsto\cX_1(t,r)$ satisfies an SDE on $[\t,t]$, and thus the It\^o's formula can be used. Moreover, for any $t\in[\t,T]$, by the second condition in (\ref{|X_1|,|X_2|}), one has
\bel{o(e)}\sup_{t\in[\t,T]}\dbE\Big[\sup_{r\in[\t,t]}|\cX_1(t,r)|^p\Big]=O(\e^{p\over2}).\ee
%
Further,
it is clear that
$$\cX_1(r,r)=X_1(r),\qq r\in[\t,T].$$
Thus,
\bel{sE(e)*}\sE(\e)=\dbE\Big[\cX_1(T,T)^\top h_{xx}\cX_1(T,T)+\int_0^T\cX_1(s,s)^\top H_{xx}(s)\cX_1(s,s)ds\Big].\ee
We now treat the terms in $\sE(\e)$. To this end, we first carry out a general calculation which will be used several times below.

\ms

Let $\Th:\D_*[\t,T]\times\Om\to\dbR^{n\times n}$ be a process such that for each $t\in[\t,T]$, $\Th(t,\cd)\in L^2_\dbF(\t,t;$ $\dbR^{n\times n})$. For each $(t,s)\in\D_*[\t,T]$, applying the martingale representation theorem to the random variable $\Th(t,s)$, one can find a unique $\L(t,s,\cd)$ such that
\bel{Pi}\Pi(t,s,r)\equiv\dbE_r[\Th(t,s)]=\Th(t,s)-\int_r^s\L(t,s,\th)dW(\th),\qq\t\les r\les  s\les t\les T.\ee
Then, applying It\^o's formula to the map $r\mapsto\Pi(t,s,r)\cX_1(s,r)$, we get
$$\ba{ll}
\ns\ds d\big[\Pi(t,s,r)\cX_1(s,r)\big]=\Big[\(\Pi(t,s,r)b_x(s,r)+\L(t,s,r)\si_x(s,r)\)
X_1(r)+\L(t,s,r)\d\si^\e(s,r)\Big]dr\\
\ns\ds\qq\qq\qq\qq\qq+\Big[\L(t,s,r)\cX_1(s,r)+\Pi(t,s,r)\(\si_x(s,r)X_1(r)+\d\si^\e(s,r)\)
\Big]dW(r).\ea$$
Further, applying the It\^o's formula to the map $r\mapsto\cX_1(t,r)^\top\Pi(t,s,r)\cX_1(s,r)$ yields
$$\ba{ll}
\ns\ds d\big[\cX_1(t,r)^\top\Pi(t,s,r)\cX_1(s,r)\big]\\
\ns\ds=\Big[X_1(r)^\top b_x(t,r)^\top\Pi(t,s,r)\cX_1(s,r)+\cX_1(t,r)^\top
\(\Pi(t,s,r)b_x(s,r)+\L(t,s,r)\si_x(s,r)\)X_1(r)\\
\ns\ds\qq+\cX_1(r,s)^\top\L(t,s,r)\d\si^\e(s,r)+\(X_1(r)^\top\si_x(t,r)^\top
+\d\si^\e(t,r)^\top\)\\
\ns\ds\qq\q\cd\(\L(t,s,r)\cX_1(s,r)+\Pi(t,s,r)\si_x(s,r)X_1(r)+\Pi(t,s,r)\d\si^\e(s,r)\)
\Big]dr+\G(t,s,r)dW(r)\\
\ns\ds=\Big[X_1(r)^\top\(b_x(t,r)^\top\Pi(t,s,r)+\si_x(t,r)^\top\L(t,s,r)\)\cX_1(s,r)\\
\ns\ds\qq+\cX_1(t,r)^\top
\(\Pi(t,s,r)b_x(s,r)+\L(t,s,r)\si_x(s,r)\)X_1(r)\\
\ns\ds\qq+X_1(r)^\top\si_x(t,r)^\top\Pi(t,s,r)\si_x(s,r)X_1(r)+\d\si^\e(t,r)^\top
\Pi(t,s,r)\d\si^\e(s,r)\\
\ns\ds\qq+\d\si^\e(t,r)^\top\(\L(t,s,r)\cX_1(s,r)+\Pi(t,s,r)\si_x(s,r)X_1(r)\)\\
\ns\ds\qq+\(\cX_1(r,s)^{\top}\L(t,s,r)+X_1(r)^\top
\si_x(t,r)^\top\Pi(t,s,r)\)\d\si^\e(s,r)\Big]dr+\G(t,s,r)dW(r),\ea$$
where
\bel{G}\ba{ll}
\ns\ds\G(t,s,r)=\(X_1(r)^\top\si_x(t,r)^\top
+\d\si^\e(t,r)^\top\)\Pi(t,s,r)\cX_1(s,r)\\
\ns\ds\qq\qq\qq+\cX_1(t,r)^\top\Big[\L(t,s,r)\cX_1(s,r)+\Pi(t,s,r)\(\si_x(s,r)X_1(r)+\d\si^\e(s,r)\)
\Big].\ea\ee
By the integrability of $(\Pi,\L)$, $X_1$, $\cX_1$, and (H1), we see that
\bel{Infinitesimal-1}\ba{ll}
\ns\ds\dbE\int_0^T\Big[\d\si^\e(t,r)^\top\(\L(t,s,r)\cX_1(s,r)+\Pi(t,s,r)\si_x(s,r)X_1(r)\)\\
\ns\ds\qq\q+\(\cX_1(r,s)^{\top}\L(t,s,r)+X_1(r)^\top
\si_x(t,r)^\top\Pi(t,s,r)\)\d\si^\e(s,r)\Big]dr=o(\e).\ea\ee
Thus, (noting $\t\les r\les s\les t\les T$, and $\cX_1(t,r)$ and $\cX_1(s,r)$ being $\cF_r$-measurable)
\bel{cXPcX}\ba{ll}
\ns\ds\dbE\big[\cX_1(t,r)^\top\Th(t,s)\cX_1(s,r)\big]=\dbE\big[\cX_1(t,r)^\top\Pi(t,s,r)\cX_1(s,r)\big]\\
\ns\ds=\dbE\int_\t^r\Big[X_1(\th)^\top\(b_x(t,\th)^\top\Pi(t,s,\th)+\si_x(t,\th)^\top
\L(t,s,\th)\)\cX_1(s,\th)\\
\ns\ds\qq+\cX_1(t,\th)^\top
\(\Pi(t,s,\th)b_x(s,\th)+\L(t,s,\th)\si_x(s,\th)\)X_1(\th)\\
\ns\ds\qq+X_1(\th)^\top\si_x(t,\th)^\top\Pi(t,s,\th)\si_x(s,\th)X_1(\th)+\d\si^\e(t,\th)^\top
\Pi(t,s,\th)\d\si^\e(s,\th)\Big]d\th+o(\e)\\
\ns\ds=\dbE\int_\t^r\Big[X_1(\th)^\top\(b_x(t,\th)^\top\Th(t,s)+\si_x(t,\th)^\top
\L(t,s,\th)\)\cX_1(s,\th)\\
\ns\ds\qq+\cX_1(t,\th)^\top
\(\Th(t,s)b_x(s,\th)+\L(t,s,\th)\si_x(s,\th)\)X_1(\th)\\
\ns\ds\qq+X_1(\th)^\top\si_x(t,\th)^\top\Th(t,s)\si_x(s,\th)X_1(\th)+\d\si^\e(t,\th)^\top
\Th(t,s)\d\si^\e(s,\th)\Big]d\th+o(\e).\ea\ee
In the above, to get rid of the It\^o's integral of $\G(t,s,\cd)$ under the expectation, some standard localization arguments via stopping times have been used. We point out that the arguments of It\^{o}'s formula on $\cX_1$ also appeared in \cite{Wang-2016-submitted}. Now, we are ready to look at the terms in $\sE(\e)$ (see \rf{sE(e)}).

\ms

\it Step 1. Treatment of the first term in $\sE(\e)$. \rm

\ms

Take $t=s=T$ and $\Th(T,T)=h_{xx}$ in \rf{Pi}. Denote the adapted solution of \rf{Pi} by
$$(\Pi(T,T,r),\L(T,T,r))\equiv(P_1(r),Q_1(r)),\qq r\in[\t,T].$$
Then \rf{Pi} reads
\bel{BSDE1}P_1(r)=h_{xx}-\int_r^T Q_1(\th)dW(\th),\qq r\in[\t,T].\ee
By the uniqueness, $ P_1(\cd)$ and $ Q_1(\cd)$ take values in $\dbS^n$.
Making use of (\ref{|X_1|,|X_2|}) and (H1), one has
$$\ba{ll}
\ns\ds\dbE\int_\t^{\t+\e} \big|X_1(s)\big| \big|P_1(s)\big|\big|\d\si(T,s)\big|ds\\
\ds\les\Big[\dbE\int_\t^{\t+\e} \big|X_1(s)\big|^4  ds\Big]^{1\over4}
\Big[\dbE\int_\t^{\t+\e} \big|\d\si(T,s)\big|^4ds\Big]^{1\over4}\Big[\dbE\int_\t^{\t+\e} \big|P_1(s)\big|^2ds\Big]^{1\over2}=o(\e),\ea$$
and
$$\ba{ll}
\ds\dbE\int_\t^{\t+\e}\big|\cX_1(T,s)\big|\big|Q_1(s)\big|\big|\d\si(T,s)\big|ds\\
\ds\les\Big[\dbE\int_\t^{\t+\e} \big|\cX_1(T,s)\big|^4ds\Big]^{1\over4}
\Big[\dbE\int_\t^{\t+\e}\big|\d\si(T,s)\big|^4  ds\Big]^{1\over4} \Big[\dbE\int_\t^{\t+\e} \big|Q_1(s)\big|^2ds\Big]^{1\over2}=o(\e).\ea$$
Therefore, the corresponding (\ref{Infinitesimal-1}) follows. On the other hand, in the current case, we have
$$\ba{ll}
\ds\G(T,T,r)=\cX_1(T,r)^{\top}Q_1(r)\cX_1(T,r)+
\cX_1(T,r)^\top P_1(r)\(\si_x(T,r)X_1(r)+\d\si^\e(T,r)\)\\
\ds \qq\qq\qq+\(X_1(r)^{\top}\si_x(T,r)^{\top}+\d\si^\e(T,r)^{\top}\)P_1(r)\cX_1(T,r).
\ea$$
Again, by a localization argument, we can get rid of the It\^o's integral of $\G(T,T,\cd)$ under expectation to obtain the following, which is a particular case of (\ref{cXPcX})
\bel{cXPcX1}\ba{ll}
\ds\dbE\big[X_1(T)^\top h_{xx}X_1(T)\big]=\dbE\big[\cX_1(T,T)^\top P_1(T)\cX_1(T,T)\big]\\
\ds=\dbE\int_\t^T\Big[X_1(r)^\top\(b_x(T,r)^\top P_1(r)+\si_x(T,r)^\top
 Q_1(r)\)\cX_1(T,r)\\
\ds\qq+\cX_1(T,r)^\top
\( P_1(r)b_x(T,r)+ Q_1(r)\si_x(T,r)\)X_1(r)\\
\ds\qq+X_1(r)^\top\si_x(T,r)^\top P_1(r)\si_x(T,r)X_1(r)+\d\si^\e(T,r)^\top
 P_1(r)\d\si^\e(T,r)\Big]dr+o(\e)\\
\ds\equiv\dbE\int_\t^T\Big[X_1(r)^\top F_1(r)\cX_1(T,r)+\cX_1(T,r)^\top
F_1(r)^\top X_1(r)\\
\ds\qq+X_1(r)^\top G_1(r)X_1(r)+\d\si^\e(T,r)^\top
 P_1(r)\d\si^\e(T,r)\Big]dr+o(\e).\ea\ee
Here,
\bel{F_1G_1}\left\{\2n\ba{ll}
\ns\ds F_1(r)=b_x(T,r)^\top P_1(r)+\si_x(T,r)^\top Q_1(r),\\
\ns\ds G_1(r)=\si_x(T,r)^\top P_1(r)\si_x(T,r).\ea\right.\ee
Note that both $F_1(r)$ and $G_1(r)$ are $\cF_r$-measurable and independent of $u^\e(\cd)$. Consequently, at the end of this step, one has
\bel{sE(e)1}\ba{ll}
\ds\sE(\e)=\dbE\int_\t^T\Big[X_1(r)^\top F_1(r)\cX_1(T,r)+\cX_1(T,r)^\top
F_1(r)^\top X_1(r)\\
\ds\qq\qq\qq+X_1(r)^\top\1n\(H_{xx}(r)\1n+\1n G_1(r)\)X_1(r)\1n+\1n\d\si^\e(T,r)^\top
\1n P_1(r)\d\si^\e(T,r)\Big]dr\1n+\1n o(\e).\ea\ee
We see that, with the above step, the term $\dbE\big[X_1(T)^\top h_{xx}X_1(T)\big]$ is absorbed. However, the new term of the form $X_1(r)^\top F_1(r)\cX_1(T,r)$ and its transpose appear under integral. This will be handled by the next step.

\ss

\it Step 2. Treatment of the term $X_1(r)^\top F_1(r)\cX_1(T,r)$ and its transpose. \rm

\ms

Take $t=T$ in \rf{Pi} and let $\Th(T,s)=\Th_2(T,r)$ be undetermined (which is $\cF_s$-measurable). By martingale representation, one has
\bel{Th=}\Th_2(T,s)=\dbE_r[\Th_2(T,s)]+\int_r^s\L_2(T,s,\th)dW(\th),\qq\t\les r\les s\les T.\ee
In this case, \rf{cXPcX} reads
\bel{cXPcX2}\ba{ll}
\ns\ds\dbE\big[\cX_1(T,r)^\top\Th_2(T,r)X_1(r)\big]\\
\ns\ds=\dbE\int_\t^r\Big[X_1(\th)^\top\(b_x(T,\th)^\top\Th_2(T,r)+\si_x(T,\th)^\top
\L_2(T,r,\th)\)\cX_1(r,\th)\\
\ns\ds\qq+\cX_1(T,\th)^\top
\(\Th_2(T,r)b_x(r,\th)+\L_2(T,r,\th)\si_x(r,\th)\)X_1(\th)\\
\ns\ds\qq+X_1(\th)^\top\si_x(T,\th)^\top\Th_2(T,r)\si_x(r,\th)X_1(\th)\1n
+\1n\d\si^\e(T,\th)^\top\Th_2(T,r)\d\si^\e(r,\th)\Big]d\th\1n+\1n o(\e)\\
\ns\ds\equiv\dbE\int_\t^r\Big[X_1(\th)^\top F_2(r,\th)\cX_1(r,\th)+\cX_1(T,\th)^\top\wt F_2(r,\th)^\top X_1(\th)\\
\ns\ds\qq+{1\over2}X_1(\th)^\top G_2(r,\th)X_1(\th)+\d\si^\e(T,\th)^\top
\Th_2(T,r)\d\si^\e(r,\th)\Big]d\th+o(\e),\ea\ee
where
\bel{F_2}\left\{\2n\ba{ll}
\ns\ds F_2(r,\th)=b_x(T,\th)^\top\Th_2(T,r)+\si_x(T,\th)^\top\L_2(T,r,\th),\\
\ns\ds\wt F_2(r,\th)=b_x(r,\th)^\top\Th_2(T,r)^\top+\si_x(r,\th)^\top \L_2(T,r,\th)^\top,\\
\ns\ds G_2(r,\th)=\si_x(T,\th)^\top\Th_2(T,r)\si_x(r,\th)+\si_x(r,\th)^\top
\Th_2(T,r)^\top\si_x(T,\th).\ea\right.\ee
Thus,
$$\ba{ll}
\ds\dbE\int_\t^T\big[\cX_1(T,r)^\top\Th_2(T,r)X_1(r)+X_1(r)^\top \Th_2(T,r)^\top\cX_1(T,r)\big]dr\\
\ds=\dbE\int_\t^T\int_\t^r\Big[X_1(\th)^\top F_2(r,\th)\cX_1(r,\th)+
\cX_1(r,\th)F_2(r,\th)^\top X_1(\th)\\
\ns\ds\qq+\cX_1(T,\th)^\top
\wt F_2(r,\th)^\top X_1(\th)+X_1(\th)^\top\wt F_2(r,\th)\cX_1(T,\th)+X_1(\th)^\top G_2(r,\th)X_1(\th)\\
\ns\ds\qq+\d\si^\e(T,\th)^\top
\Th_2(T,r)\d\si^\e(r,\th)+\d\si^\e(r,\th)^\top
\Th_2(T,r)^\top\d\si^\e(T,\th) \Big]d\th dr+o(\e)\\
\ds=\dbE\int_\t^T\int_r^T\Big[X_1(r)^\top F_2(\th,r)\cX_1(\th,r)+
\cX_1(\th,r)F_2(\th,r)^\top X_1(r)\\
\ns\ds\qq+\cX_1(T,r)^\top
\wt F_2(\th,r)^\top X_1(r)+X_1(r)^\top\wt F_2(\th,r)\cX_1(T,r)+X_1(r)^\top G_2(\th,r)X_1(r)\\
\ns\ds\qq+\d\si^\e(T,r)^\top
\Th_2(T,\th)\d\si^\e(\th,r)+\d\si^\e(\th,r)^\top
\Th_2(T,\th)^\top\d\si^\e(T,r)\Big]d\th dr+o(\e).\ea$$
Consequently, \rf{sE(e)1} becomes
%
$$\ba{ll}
\ns\ds\sE(\e)\equiv\dbE\Big[X_1(T)^\top h_{xx}X_1(T)+\int_\t^TX_1(r)^\top H_{xx}(r)X_1(r)dr\Big]\\
\ns\ds\qq=\dbE\int_\t^T\Big[X_1(r)^\top\(F_1(r)-\Th_2(T,r)^\top+\int_r^T\wt F_2(\th,r)d\th\)\cX_1(T,r)\\
\ns\ds\qq\q+\cX_1(T,r)^\top\(F_1(r)^\top-\Th_2(T,r)+\int_r^T\wt F_2(\th,r)^\top d\th\)X_1(r)\\
\ns\ds\qq\q+\int_r^T \big[X_1(r)^\top F_2(\th,r)\cX_1(\th,r)+
\cX_1(\th,r)^\top F_2(\th,r)^\top X_1(r)\big]d\th\\
\ns\ds\qq\q+X_1(r)^\top\(H_{xx}(r)+G_1(r)+\int_r^TG_2(\th,r)d\th\)X_1(r)+\d\si^\e(T,r)^\top
 P_1(r)\d\si^\e(T,r)\\
\ns\ds\qq\q+\2n\int_r^T\3n\big[\d\si^\e(T,r)^\top
\Th_2(T,\th)\d\si^\e(\th,r)\1n+\1n\d\si^\e(\th,r)^\top
\Th_2(T,\th)^\top\d\si^\e(T,r)\big]d\th\Big]dr\1n+\1n o(\e).\ea$$
Thus, in order to eliminate the term of the form $X_1(r)^\top(\cds)\cX_1(T,r)$ and its transpose, we let
$$\ba{ll}
\ds\Th_2(T,r)=F_1(r)^\top+\dbE_r\int_r^T\wt F_2(\th,r)^\top d\th\\
\ds\qq\q=\dbE_r\Big[P_1(r)b_x(T,r)+Q_1(r)\si_x(T,r)+\2n\int_r^T\3n\(\Th_2(T,\th)b_x(\th,r)+
\L_2(T,\th,r)\si_x(\th,r)\)d\th\Big].\ea$$
Now, by denoting
$$P_2(r)=\Th_2(T,r)^\top,\qq Q_2(\th,r)=\L_2(T,\th,r)^\top,\qq\t\les r\les\th\les T,$$
we have
$$\ba{ll}
\ns\ds P_2(r)=\dbE_r\Big[b_x(T,r)^\top P_1(r)+\si_x(T,r)^\top Q_1(r) +\int_r^T\(b_x(\th,r)^\top P_2(\th)+\si_x(\th,r)^\top Q_2(\th,r)\)d\th\Big]\\
\ns\ds\qq=b_x(T,r)^\top P_1(r)+\si_x(T,r)^\top Q_1(r)+\int_r^T\(b_x(\th,r)^\top P_2(\th)+\si_x(\th,r)^\top Q_2(\th,r)\)d\th\\
\ns\ds\qq\qq-\int_r^T\wt Q_2(r,\th)dW(\th),\ea$$
for some $\wt Q_2(r,\th)$ defined on $\D^*[\t,T]$, or for $\t\les r\les\th\les T$, with a suitable measurability and integrability. At the same time, \rf{Th=} reads
$$P_2(s)=\dbE_r[P_2(s)]+\int_r^sQ_2(s,\th)dW(\th),\qq\t\les r\les s\les T,$$
with $Q_2(s,\th)$ defined on $\D_*[t,T]$, or for $\t\les\th\les s\les T$. We now extend $Q_2(s,\th)$ from $\D_*[\t,T]$ to $[\t,T]^2$ by
$$Q_2(s,\th)=\wt Q_2(s,\th),\qq\t\les s\les\th\les T.$$
Then $(P_2(\cd),Q_2(\cd\,,\cd))$ is the unique adapted M-solution to the following Type-II BSVIE:
\bel{BSVIE2}\ba{ll}
\ns\ds P_2(r)\1n=\1n b_x(T,r)^\top P_1(r)\1n+\1n\si_x(T,r)^\top Q_1(r)\1n+\2n\int_r^T\3n\(b_x(\th,r)^\top P_2(\th)\1n+\1n\si_x(\th,r)^\top Q_2(\th,r)\)d\th\\
\ns\ds\qq\qq\qq-\int_r^TQ_2(r,\th)dW(\th),\qq r\in[\t,T],\ea\ee
which is independent of spike variation $u^\e(\cd)$ of $\bar u(\cd)$. Having the above result, we obtain:
\bel{sE**}\ba{ll}
\ds\sE(\e) =\1n\dbE\2n\int_\t^T\2n\Big[X_1(r)^\top\(H_{xx}(r)\1n+\1n G_1(r)\1n+\1n\int_r^T\3n G_2(\th,r)d\th\)X_1(r)\\
\ds\qq\qq+\int_r^T\(X_1(r)^\top F_2(\th,r)\cX_1(\th,r)+
\cX_1(\th,r)^\top F_2(\th,r)^\top X_1(r)\)d\th\Big]dr\\
\ds\qq\qq+\dbE\int_\t^T\Big[\d\si^\e(T,r)^\top P_1(r)\d\si^\e(T,r)\\
\ns\ds\qq\qq+\int_r^T\3n\(\d\si^\e(T,r)^\top
P_2(\th)^\top\d\si^\e(\th,r)+\d\si^\e(\th,r)^\top
P_2(\th)\d\si^\e(T,r)\)d\th\Big]dr+o(\e).\ea\ee
We see by \rf{F_2} that
\bel{F_2G_2}\left\{\2n\ba{ll}
\ds F_2(r,\th)=b_x(T,\th)^\top P_2(r)^\top+\si_x(T,\th)^\top Q_2(r,\th),\\
\ns\ds G_2(r,\th)=\si_x(T,\th)^\top P_2(r)^\top\si_x(r,\th)+\si_x(r,\th)^\top
P_2(r)\si_x(T,\th).\ea\right.\ee
Now, in $\sE(\e)$ (see \rf{sE**}), the terms of the forms $X_1(r)^\top(\cds)X_1(r)$ and $X_1(r)^\top(\cds)\cX_1(\th,r)$, together with its transpose, need to be handled. In the next step, we will get rid of those terms simultaneously.

\ss

\it Step 3. Treatment of the terms $X_1(r)^\top(\cds)X_1(r)$, $\cX_1(\th,r)^\top(\cds)X_1(r)$, and its transpose. \rm

\ms

To treat $X_1(r)^\top(\cds)X_1(r)$, we take $t=s$ in \rf{Pi} with $\Th(s,s)=\Th_3(s)$ taking values in $\dbS^n$, undetermined. Again, by martingale representation theorem, we have a unique $\L_3(s,s,\cd)$ such that
\bel{BSDE3}\Th_3(s)=\dbE_r[\Th_3(s)]+\int_r^s\L_3(s,s,\th)dW(\th),\qq\t\les r\les s\les T.\ee
Since $\Th_3(s)$ is symmetric, so is $\L_3(s,s,r)$. Now, \rf{cXPcX} becomes
$$\ba{ll}
\ds\dbE\big[X_1(r)^\top \Th_3(r)X_1(r)\big]=\dbE\int_\t^r\Big[X_1(\th)^\top\(b_x(r,\th)^\top\Th_3(r)+\si_x(r,\th)^\top
\L_3(r,r,\th)\)\cX_1(r,\th)\\
\ds\qq+\cX_1(r,\th)^\top
\(\Th_3(r)b_x(r,\th)+\L_3(r,r,\th)\si_x(r,\th)\)X_1(\th)\\
\ds\qq+X_1(\th)^\top\si_x(r,\th)^\top\Th_3(r)\si_x(r,\th)X_1(\th)
+\d\si^\e(r,\th)^\top
\Th_3(r)\d\si^\e(r,\th)\Big]d\th+o(\e)\\
\ds=\dbE\int_\t^r\Big[X_1(\th)^\top F_3(r,\th)\cX_1(r,\th)+\cX_1(r,\th)^\top
F_3(r,\th)^\top X_1(\th)\\
\ds\qq+X_1(\th)^\top G_3(r,\th)X_1(\th)+\d\si^\e(r,\th)^\top
\Th_3(r)\d\si^\e(r,\th)\Big]d\th+o(\e),\ea$$
where
\bel{F_3G_3}\left\{\2n\ba{ll}
\ds F_3(r,\th)=b_x(r,\th)^\top\Th_3(r)+\si_x(r,\th)^\top\L_3(r,r,\th),\\
\ds G_3(r,\th)=\si_x(r,\th)^\top\Th_3(r)\si_x(r,\th).\ea\right.\ee
Then integrate it over $[\t,T]$, one has
\bel{XThX*}\ba{ll}
\ds\dbE\2n\int_\t^T\3n\big[X_1(r)^{\1n\top}\1n\Th_3(r)X_1(r)\big]dr\1n=\1n\dbE\2n\int_\t^T\3n
\int_\t^r\3n\Big[X_1(\th)^{\1n\top}\1n F_3(r,\th)\cX_1(r,\th)\1n+\1n\cX_1(r,\th)^{\1n\top}\1n
F_3(r,\th)^{\1n\top}\1n X_1(\th)\\
\ds\qq+X_1(\th)^\top G_3(r,\th)X_1(\th)+\d\si^\e(r,\th)^\top
\Th_3(r)\d\si^\e(r,\th)\Big]d\th dr+o(\e)\\
\ds=\dbE\int_\t^T\int_r^T\Big[X_1(r)^\top F_3(\th,r)\cX_1(\th,r)+\cX_1(\th,r)^\top F_3(\th,r)^\top X_1(r)\\
\ds\qq+X_1(r)^\top G_3(\th,r)X_1(r)+\d\si^\e(\th,r)^\top
\Th_3(\th)\d\si^\e(\th,r)\Big]d\th dr+o(\e).\ea\ee
Next, for the term $\cX_1(\th,r)^\top(\cds)X_1(r)$, we let $\Th(t,s)=\Th_4(t,s)$ be undetermined and by martingale representation theorem,
$$\Th_4(t,s)=\dbE_r[\Th_4(t,s)]+\int_r^s\L_4(t,s,\th)dW(\th),\qq\t\les r\les s\les t\les T.$$
Then, \rf{cXPcX} reads
$$\ba{ll}
\ds\dbE\big[\cX_1(\th,r)^\top \Th_4(\th,r)X_1(r)\big]\\
\ds=\dbE\int_\t^r\Big[X_1(\th')^\top\(b_x(\th,\th')^\top\Th_4(\th,r)
+\si_x(\th,\th')^\top\L_4(\th,r,\th')\)\cX_1(r,\th')\\
\ds\qq+\cX_1(\th,\th')^\top
\(\Th_4(\th,r)b_x(r,\th')+\L_4(\th,r,\th')\si_x(r,\th')\)X_1(\th')\\
\ds\qq+X_1(\th')^\top\si_x(\th,\th')^\top\Th_4(\th,r)\si_x(r,\th')X_1(\th')
+\d\si^\e(\th,\th')^\top\Th_4(\th,r)\d\si^\e(r,\th')\Big]d\th'+o(\e)\\
\ds\equiv\dbE\int_\t^r\Big[X_1(\th')^\top F_4(\th,r,\th')\cX_1(r,\th')+\cX_1(\th,\th')^\top\wt F_4(\th,r,\th')^\top X_1(\th')\\
\ds\qq+{1\over2}X_1(\th')^\top G_4(\th,r,\th')X_1(\th')
+\d\si^\e(\th,\th')^\top\Th_4(\th,r)\d\si^\e(r,\th')\Big]d\th'+o(\e),\ea$$
where (note that $\t\les\th'\les r\les\th\les T$)
\bel{F_4G_4}\left\{\2n\ba{ll}
\ns\ds F_4(\th,r,\th')=b_x(\th,\th')^\top\Th_4(\th,r)+\si_x(\th,\th')^\top \L_4(\th,r,\th'),\\
\ns\ds\wt F_4(\th,r,\th')=b_x(r,\th')^\top\Th_4(\th,r)^\top+\si_x(r,\th')^\top\L_4(\th,r,\th')^\top,\\
\ns\ds G_4(\th,r,\th')=\si_x(\th,\th')^\top\Th_4(\th,r)\si_x(r,\th')
+\si_x(r,\th')^\top\Th_4(\th,r)^\top\si_x(\th,\th').\ea\right.\ee
Therefore, we have:
$$\ba{ll}
\ns\ds\dbE\int_\t^T\int_r^T\cX_1(\th,r)^\top\Th_4(\th,r)X_1(r)d\th dr\\
\ns\ds=\dbE\int_\t^T\int_r^T\int_\t^r\Big[X_1(\th')^\top  F_4(\th,r,\th')\cX_1(r,\th')+\cX_1(\th,\th')^\top\wt F_4(\th,r,\th')^\top X_1(\th')\\
\ns\ds\qq+{1\over2}X_1(\th')^\top G_4(\th,r,\th')X_1(\th')
+\d\si^\e(\th,\th')^\top\Th_4(\th,r)\d\si^\e(r,\th')\Big]d\th'd\th dr+o(\e)\\
\ns\ds=\dbE\int_\t^T\int_{\th'}^T\int_\t^{\th'}\Big[X_1(r)^\top  F_4(\th,\th',r)\cX_1(\th',r)+\cX_1(\th,r)^\top\wt F_4(\th,\th',r)^\top X_1(r)\\
\ns\ds\qq+{1\over2}X_1(r)^\top G_4(\th,\th',r)X_1(r)
+\d\si^\e(\th,r)^\top\Th_4(\th,\th')\d\si^\e(\th',r)\Big]drd\th d\th'+o(\e)\\
\ns\ds=\dbE\int_\t^T\int_\t^{\th'}\int_{\th'}^T\Big[X_1(r)^\top  F_4(\th,\th',r)\cX_1(\th',r)+\cX_1(\th,r)^\top\wt F_4(\th,\th',r)^\top X_1(r)\\
\ns\ds\qq+{1\over2}X_1(r)^\top G_4(\th,\th',r)X_1(r)
+\d\si^\e(\th,r)^\top\Th_4(\th,\th')\d\si^\e(\th',r)\Big]d\th drd\th'+o(\e)\\
\ns\ds=\dbE\int_\t^T\int_r^T\int_{\th'}^T\Big[X_1(r)^\top  F_4(\th,\th',r)\cX_1(\th',r)+\cX_1(\th,r)^\top\wt F_4(\th,\th',r)^\top X_1(r)\\
\ns\ds\qq+{1\over2}X_1(r)^\top G_4(\th,\th',r)X_1(r)
+\d\si^\e(\th,r)^\top\Th_4(\th,\th')\d\si^\e(\th',r)\Big]d\th d\th'dr+o(\e)\\
\ns\ds=\1n\dbE\2n\int_\t^T\3n\Big[\1n\int_r^T\3n\int_\th^T\3n X_1(r)^{\1n\top}\1n  F_4(\th',\th,r)\cX_1(\th,r)d\th'd\th\1n+\2n\int_r^T\3n\int_r^\th\3n
\cX_1(\th,r)^\top\wt F_4(\th,\th',r)^\top X_1(r)d\th'd\th\\
\ns\ds\qq+\2n\int_r^T\3n\int_\th^T\3n\({1\over2}X_1(r)^{\1n\top}\1n G_4(\th',\th,r)X_1(r)
\1n+\1n\d\si^\e(\th',r)^{\1n\top}\1n\Th_4(\th',\th)\d\si^\e(\th,r)\)d\th'd\th\Big]dr\1n
+\1no(\e).\ea$$
Hence, combining \rf{sE**} and \rf{XThX*} with the above, we have
%
$$\ba{ll}
\ns\ds\sE(\e)\\
\ns\ds=\1n\dbE\2n\int_\t^T\3n\Big\{X_1(r)^{\1n\top}\1n\Big[H_{xx}(r)\1n+\1n G_1(r)\1n+\2n\int_r^T\3n\(G_2(\th,r)\1n+\1n G_3(\th,r)\1n+\2n\int_\th^T\3n G_4(\th',\th,r)d\th'\)d\th\1n-\1n\Th_3(r)\Big]X_1(r)\\
\ns\ds\q+\2n\int_r^T\3n\Big[X_1(r)^{\1n\top}\1n\(F_2(\th,r)\1n+\1n F_3(\th,r)\1n-\1n\Th_4(\th,r)^{\1n\top}\2n+\2n\int_\th^T
\3n F_4(\th',\th,r)d\th'\1n+\2n\int_r^\th\3n\wt F_4(\th,\th',r)d\th'\)\cX_1(\th,r)\\
\ns\ds\q+\cX_1(\th,r)^{\1n\top}\1n\(F_2(\th,r)\1n+\1n F_3(\th,r)\1n-\1n\Th_4(\th,r)^{\1n\top}\2n+\2n\int_\th^T
\3n F_4(\th',\th,r)d\th'\1n+\2n\int_r^\th\3n\wt F_4(\th,\th',r)d\th'\)^{\2n\top}\2n X_1(r)\Big]d\th\\
\ns\ds\q+\d\si^\e(T,r)^\top P_1(r)\d\si^\e(T,r)\\
\ns\ds\q+\int_r^T\3n\Big[\d\si^\e(T,r)^\top
P_2(\th)^\top\d\si^\e(\th,r)+\d\si^\e(\th,r)^\top
P_2(\th)\d\si^\e(T,r)+\d\si^\e(\th,r)^\top
\Th_3(\th)\d\si^\e(\th,r)\\
\ns\ds\q+\int_\th^T\3n\(\d\si^\e(\th',r)^\top\Th_4(\th',\th)\d\si^\e(\th,r)+\d\si^\e(\th,r)^\top
\Th_4(\th',\th)^\top
\d\si^\e(\th',r)\)d\th'\Big]d\th\Big\}dr\1n+\1n o(\e).\ea$$
Consequently, in order to eliminate the term $X_1(r)^\top(\cds)X_1(r)$, we should let
$$\ba{ll}
\ds \Th_3(r)\1n=\1nH_{xx}(r)\1n+\1n G_1(r)\1n+\2n\int_r^T\3n\(G_2(\th,r)\1n+\1n G_3(\th,r)\1n+\2n\int_\th^T\3n G_4(\th',\th,r)d\th'\)d\th\1n-\2n\int_r^T\3n\wt Q_3(r,\th)dW(\th)\\
\ds\qq=H_{xx}(r)+\si_x(T,r)^\top P_1(r)\si_x(T,r)+\int_r^T\3n\(\si_x(T,r)^\top P_2(\th)^\top\si_x(\th,r)\\
\ds\qq\qq+\si_x(\th,r)^\top
P_2(\th)\si_x(T,r)+\si_x(\th,r)^\top\Th_3(\th)\si_x(\th,r)\)d\th\\
\ns\ds\qq\qq+\int_r^T\3n\int_\th^T\3n\(\si_x(\th',r)^\top\Th_4(\th',\th)\si_x(\th,r)
+\si_x(\th,r)^\top\Th_4(\th',\th)^\top\si_x(\th',r)\)d\th'd\th\\
 \ns\ds\qq\qq-\int_r^T\wt Q_3(r,\th)dW(\th),\q\t\les r\les T,\ea$$
for some $\wt Q_3(r,\th)$ defined for $\t\les r\les\th\les T$. Similarly, to eliminate $X_1(r)^\top(\cds)\cX_1(\th,r)$ and its transpose, we should let
$$\ba{ll}
\ns\ds\Th_4(\th,r)^\top\2n=\1n F_2(\th,r)\1n+\1n F_3(\th,r)+\2n\int_\th^T
\3n F_4(\th',\th,r)d\th'\1n+\2n\int_r^\th\3n\wt F_4(\th,\th',r)d\th'\1n-\2n\int_r^T\3n\wt Q_4(\th,r,\th')dW(\th')\\
\ns\ds=b_x(T,r)^\top P_2(\th)^\top +\si_x(T,r)^\top Q_2(\th,r)^\top+b_x(\th,r)^\top\Th_3(\th) +\si_x(\th,r)^\top\L_3(\th,\th,r)\\
\ns\ds\qq+\int_\th^T\(b_x(\th',r)^\top \Th_4(\th',\th)+\si_x(\th',r)^\top\L_4(\th',\th,r)\)d\th'\\
\ns\ds\qq+\int_r^\th\(b_x(\th',r)^\top\Th_4(\th,\th')^\top+\si_x(\th',r)^\top\L_4(\th,\th',r)^\top\)d\th'-\int_r^T\wt Q_4(\th,r,\th')dW(\th'),\\
\ns\ds\qq\qq\qq\qq\qq\qq\qq\qq\qq\qq\t\les r\les\th\les T,\ea$$
for some $\wt Q_4(\th,r,\th')$ defined for $\t\les r\les\th'\les T$ and $\t\les r\les\th\les T$. Let
$$\ba{ll}
\ns\ds P_3(\th)=\Th_3(\th),\q Q_3(\th,r)=\left\{\2n\ba{ll}
\ns\ds\L_3(\th,\th,r),\q\t\les r\les\th\les T,\\
\ns\ds\wt Q_3(\th,r),\q\t\les \th\les r\les T,\ea\right.\\
\ns\ds P_4(\th,\th')=\Th_4(\th,\th')^{\top},\q Q_4(\th,r,\th')=\left\{\2n\ba{ll}
\ns\ds\L_4(\th,r,\th')^{\top},\qq\t\les\th'\les r\les\th\les T,\\
\ns\ds\wt Q_4(\th,r,\th')^\top,\qq\t\les r\les\th'\les T,\ea\right.\ea$$
Then,
\bel{P_3}\ba{ll}
\ns\ds P_3(r)=H_{xx}(r)+\si_x(T,r)^\top P_1(r)\si_x(T,r)+\int_r^T\3n\(\si_x(T,r)^\top P_2(\th)^\top\si_x(\th,r)\\
\ns\ds\qq\qq+\si_x(\th,r)^\top
P_2(\th)\si_x(T,r)+\si_x(\th,r)^\top P_3(\th)\si_x(\th,r)\)d\th\\
\ns\ds\qq\qq+\int_r^T\3n\int_\th^T\3n\(\si_x(\th',r)^\top
P_4(\th',\th)^\top\si_x(\th,r)
+\si_x(\th,r)^\top P_4(\th',\th)\si_x(\th',r)\)d\th'd\th\\
\ns\ds\qq\qq-\int_r^TQ_3(r,\th)dW(\th),\qq
\t\les r\les T,\ea\ee
and
\bel{P_4}\ba{ll}
\ns\ds P_4(\th,r)=b_x(T,r)^\top P_2(\th)^\top\2n+\si_x(T,r)^\top Q_2(\th,r)^\top\2n+b_x(\th,r)^\top P_3(\th)+\si_x(\th,r)^\top Q_3(\th,r)\\
\ns\ds\qq\qq+\2n\int_\th^T\3n\(b_x(\th',r)^\top P_4(\th',\th)^\top+\si_x(\th',r)^\top
Q_4(\th',\th,r)^\top\)d\th'\\
\ns\ds\qq\qq+\2n\int_r^\th\3n\(b_x(\th',r)^\top P_4(\th,\th')+\si_x(\th',r)^\top Q_4(\th,\th',r)\)d\th'\2n-\2n\int_r^T\3n Q_4(\th,r,\th')dW(\th'),\\
\ns\ds\qq\qq\qq\qq\qq\qq\qq\qq\qq\qq\t\les r\les\th\les T.\ea\ee
Thus, we obtain the following system
\bel{Second-order adjoint equations}\left\{\2n\ba{ll}
\ns\ds P_1(t)=h_{xx}-\int_t^T Q_1(s)dW(s),\qq\t\les t\les T,\\
\ns\ds P_2(t)\1n=\1n b_x(T,t)^{\1n\top}\1n P_1(t)\1n+\1n\si_x(T,t)^{\1n\top}\1n Q_1(t)\1n+\2n\int_t^T\3n\(b_x(s,t)^{\1n\top}\1n P_2(s)\1n+\1n\si_x(s,t)^{\1n\top}\1n Q_2(s,t)\)ds\\
\ns\ds\qq\qq\qq-\int_t^TQ_2(t,s)dW(s),\qq\t\les t\les T,\\
%
%
%
%
\ns\ds P_3(t)=H_{xx}(t)+\si_x(T,t)^\top P_1(t)\si_x(T,t)\\
\ns\ds\qq\qq+\int_t^T\3n\Big[\si_x(T,t)^\top P_2(s)^\top\si_x(s,t)+\si_x(s,t)^\top
P_2(s)\si_x(T,t)\\
\ns\ds\qq\qq\qq+\int_t^T\3n\(\si_x(\th,t)^\top P_4(\th,s)^\top\si_x(s,t)+\si_x(s,t)^\top\1n P_4(\th,s)\si_x(\th,t)\)d\th\Big]ds\\
\ns\ds\qq\qq+\int_t^T\3n\si_x(s,t)^\top P_3(s)\si_x(s,t)ds-\2n\int_t^TQ_3(t,s)dW(s),\q
\t\les t\les T\\
%
%
\ns\ds P_4(r,t)\1n=\1n b_x(T,t)^\top\1n P_2(r)^\top\3n+\1n\si_x(T,t)^\top\1n Q_2(r,t)^\top\3n+\1n b_x(r,t)^\top\1n P_3(r)\1n+\1n\si_x(r,t)^\top Q_3(r,t)\\
\ns\ds\qq\qq\qq+\int_r^T\(b_x(s,t)^\top P_4(s,r)^\top+\si_x(s,t)^\top Q_4(s,r,t)^\top\)ds\\
\ns\ds\qq\qq\qq+\int_t^r\(b_x(s,t)^\top P_4(r,s)+\si_x(s,t)^\top Q_4(r,s,t)\)ds\\
\ns\ds\qq\qq\qq-\int_t^TQ_4(r,t,s)dW(s),\qq\t\les t\les r\les T.\ea\right.\ee
By saying $(P_i,Q_i)$ to be the adapted M-solutions of the corresponding BSVIEs, the following constraints must hold:
\bel{Second-order-adjoint-constraints}\left\{\2n\ba{ll}
\ns\ds P_2(t)=\dbE_\th[P_2(t)]+\int_\th^tQ_2(t,s)dW(s),\qq\t\les\th\les t\les T,\\
\ns\ds P_3(t)=\dbE_\th[P_3(t)]+\int_\th^tQ_3(t,s)dW(s),\qq\t\les\th\les t\les T,\\
\ns\ds P_4(r,t)=\dbE_\th[P_4(r,t)]+\int_\th^tQ_4(r,t,s)dW(s),\qq\t\les \th\les t\les r\les T.\ea\right.\ee
The above is called the {\it second order adjoint equation}. With the above system, we eventually end up with the following representation theorem.

\bt\label{variation-2} \sl Let {\rm(H1)--(H2)} hold and $(\bar X(\cd),\bar u(\cd))$ be a given state-control pair. Let \rf{Second-order adjoint equations} admit an adapted M-solution $(P_i,Q_i)$ ($1\les i\les4$) such that
$$\ba{ll}
\ns\ds (P_1,Q_1)\in L^2_{\dbF}(\Omega;C([\t,T];\dbS^n))\times L^2_{\dbF}(\t,T;\dbS^n),\\
\ns\ds  (P_2,P_3,P_4)\in  L^2_{\dbF}(\t,T;\dbR^{n\times n})\times L^2_{\dbF}(\t,T;\dbS^{n})\times L^2_\dbF(\Om;L^2(\D_*[\t,T];\dbR^{n\times n})).\ea$$
Then $\sE(\e)$ admits the following representation.
$$\ba{ll}
\ns\ds\sE(\e)=\1n\dbE\2n\int_\t^T\3n\Big\{\d\si^\e(T,r)^\top P_1(r)\d\si^\e(T,r)+\int_r^T\3n\Big[\d\si^\e(T,r)^\top
P_2(\th)^\top\d\si^\e(\th,r)\\
\ns\ds\qq\qq\q+\d\si^\e(\th,r)^\top
P_2(\th)\d\si^\e(T,r)+\d\si^\e(\th,r)^\top P_3(\th)\d\si^\e(\th,r)\\
\ns\ds\qq\qq\q+\int_\th^T\3n\(\d\si^\e(\th',r)^\top P_4(\th',\th)^\top\d\si^\e(\th,r)+\d\si^\e(\th,r)^\top
P_4(\th',\th)\d\si^\e(\th',r)\)d\th'\Big]d\th\Big\}dr\\
\ns\ds\qq\qq\q+o(\e).\ea$$

\et

This will lead to the maximum principle and so on, whose precise statements will be given in Section 6.

\section{Well-posedness of the Second Order Adjoint Equations}

In this section, we will establish the well-posedness of system (\ref{Second-order adjoint equations}). To begin with, let us make some observations. The existence and uniqueness of $(P_1,Q_1)$ is easy to see, while the existence and uniqueness of $(P_2,Q_2)$ follows from the BSVIEs theory in \cite{Yong-2008}.
The remaining solvability of (\ref{Second-order adjoint equations}) lies in $(P_3,Q_3,P_4,Q_4)$.
For the equation of $(P_3,Q_3)$, it is a classical linear BSVIE with given $(P_1,P_2,P_4)$. However, as to the equation of $(P_4,Q_4)$, its well-posedness cannot be directly given by the current BSVIEs theory. Therefore, we need to carefully treat it.

\ms

Note that by applying conditional expectation operator $\dbE_t$ on both sides of fourth equation in (\ref{Second-order adjoint equations}), we end up with the following form system (noting the third equality in \rf{Second-order-adjoint-constraints})
\bel{Second-order-1}\left\{\2n\ba{ll}
\ns\ds P(r,t)=F(r,t)+\int_r^T\big[
b_x(s,t)^\top \dbE_t P(s,r)^{\top}+\si_x(s,t)^\top Q(s,r,t)^{\top}\big]ds\\
\ns\ds\qq\qq\qq+\int_t^r\big[b_x(s,t)^\top \dbE_t P(r,s)+\si_x(s,t)^\top Q(r,s,t)\big]ds,\q r\ges t,\\
\ns\ds \dbE_\th P(r,t)=P(r,t)-\int_\th^{t} Q(r,t,s)dW(s),\  \ \t\les \th\les t\les r\les T,\ea\right.\ee
with $P(r,t)=P_4(r,t)$, $Q(r,t,s)=Q_4(r,t,s)$ and
\bel{F}F(r,t)=\dbE_t\Big[b_x(T,t)^\top\1n P_2(r)^\top\3n+\1n\si_x(T,t)^\top\1n Q_2(r,t)^\top\3n+\1n b_x(r,t)^\top\1n P_3(r)\1n+\1n\si_x(r,t)^\top Q_3(r,t)\Big].\ee
The following gives the wellposedness for (\ref{Second-order-1}).
\begin{lemma}\label{Preparation-1}
\sl Let $b_x(\cd\,,\cd),\si_x(\cd\,,\cd)$ be bounded and $F(\cd\,,\cd)\in L^2_\dbF(\Om;L^2(\D_*[\t,T]))$. Then  \rf{Second-order-1} admits a unique solution $P(\cd\,,\cd)\in  L^2_\dbF(\Om;L^2(\D_*[\t,T]))$. In addition, for any $\rho\in[\t,T]$, the following estimate holds:
\bel{Estimate-1-preparation}\ba{ll}
\ns\ds \dbE\int_{\rho}^T \int_{\rho}^r  | P(r,t)|^2 dtdr \les K \dbE\int_{\rho}^T \int_{\rho}^r | F(r,t)|^2 dtdr,
\ea\ee
where $K$ only depends on $\|b_x\|_{\infty}$, $\|\si_x\|_{\infty}$ and $T$, but not on $\rho$.
\end{lemma}

\it Proof. \rm Let $\rho\in[\t,T]$ be fixed. For any $p(\cd\,,\cd)\in  L^2_\dbF(\Om;L^2(\D_*[\rho,T]))$, we look at the following system
\bel{Second-order-1-contraction}\left\{\ba{ll}
\ns\ds P(r,t)=F(r,t)+\int_r^T\big[
b_x(s,t)^\top \dbE_t p(s,r)^{\top}+\si_x(s,t)^\top q(s,r,t)^{\top}\big]ds\\
\ns\ds\qq\qq+\int_t^r\big[b_x(s,t)^\top \dbE_t p(r,s)+\si_x(s,t)^\top q(r,s,t)\big]ds,\ \ r\ges t,\\
\ns\ds \dbE_\th p(r,t)=p(r,t)-\int_\th^{t} q(r,t,s)dW(s),\qq\rho\les \th \les t\les r\les T.\ea\right.\ee
By the first equality in (\ref{Second-order-1-contraction}), we have (for any $\b>0$)
\bel{estimate}\ba{ll}
\ns\ds \dbE\int_{\rho}^T \int_{\rho}^r e^{\b(r+t)}| P(r,t)|^2 dt dr\les  K \dbE\int_{\rho}^T \int_{\rho}^r e^{\b(r+t)}| F(r,t)|^2 dt dr\\
\ns\ds\q+K\big\|b_x(\cd\,,\cd)\big\|^2_{\infty}\dbE\int_{\rho}^T \int_{\rho}^r  e^{\b(r+t)}\Big\{\(\int_t^r \big|p(r,s)\big|ds\)^2+\(\int_r^T\big|p(s,r)\big|ds\)^2\Big\}dt dr\\
\ns\ds\q+K \big\|\si_x(\cd\,,\cd)\big\|^2_{\infty}\dbE\int_{\rho}^T \int_{\rho}^r  e^{\b(r+t)}\Big\{\(\int_t^r \big|q(r,s,t)\big|ds\)^2
+\(\int_r^T \big|q(s,r,t)\big|ds\)^2\Big\}dt dr.\ea\ee
We now estimate the right-hand side of \rf{estimate}. For the second term, by Fubini theorem, we have
%
$$\ba{ll}
\ns\ds\dbE\int_\rho^T\int_\rho^re^{\b(r+t)}\Big\{\(\int_t^r\big|  p(r,s)\big|ds\)^2+\(\int_r^T\big|p(s,r)\big|ds\)^2\Big\}dtdr\\
\ns\ds \les \dbE\int_\rho^T\int_\rho^re^{\b(r+t)}\Big\{\(\int_t^r e^{-\b s}ds\)\(\int_t^r e^{\b s}\big|p(r,s)\big|^2ds\)\\
\ns\ds\qq\qq\qq\qq\q+\(\int_r^T e^{-\b s}ds\)\(\int_r^Te^{\b s}\big|p(s,r)\big|^2ds\)\Big\}dt dr\\
\ns\ds=\dbE\int_\rho^T\int_\rho^r\int_t^re^{\b(r+t)}{e^{-\b t}-e^{-\b r}\over\b}e^{\b s}|p(r,s)|^2dsdtdr\\
\ns\ds\qq+\dbE\int_\rho^T\int_\rho^r\int_r^Te^{\b(r+t)}{e^{-\b r}-e^{-\b T}\over\b}e^{\b s}|p(s,r)|^2dsdtdr\ea$$
$$\ba{ll}
\ns\ds={1\over\b}\dbE\int_\rho^T\int_\rho^r\int_\rho^s\(e^{\b(r+s)}-e^{\b(t+s)}\)|p(r,s)|^2dtdsdr\\
\ns\ds\qq+{1\over\b}\dbE\int_\rho^T\int_r^T\int_\rho^r\(e^{\b(t+s)}-e^{\b(r+s+t-T)}\)|p(s,r)|^2dtdsdr\\
\ns\ds={1\over\b}\dbE\int_\rho^T\int_\rho^r\int_\rho^s\(e^{\b(r+s)}-e^{\b(t+s)}\)|p(r,s)|^2dtdsdr\\
\ns\ds\qq+{1\over\b}\dbE\int_\rho^T\int_\rho^s\int_\rho^r\(e^{\b(t+s)}-e^{\b(r+s+t-T)}\)|p(s,r)|^2
dtdrds\\
\ns\ds\les{2T\over\b}\dbE\int_\rho^T \int_\rho^r  e^{\b (r+s)}|p(r,s)|^2dsdr.\ea$$
For the third term on the right-hand side of \rf{estimate}, similarly, by Fubini theorem,
%
$$\ba{ll}
\ns\ds\dbE\int_\rho^T\int_\rho^re^{\b(r+t)}\Big\{\(\int_t^r\big|q(r,s,t)\big|ds\)^2
+\(\int_r^T \big|q(s,r,t)\big|ds\)^2\Big\}dtdr\\
\ns\ds\les\dbE\int_\rho^T\int_\rho^re^{\b(r+t)}\Big\{\(\int_t^r e^{-\b s}ds\)\(\int_t^r e^{\b s}\big| q(r,s,t)\big|^2ds\)\\
\ns\ds\qq\qq\qq\qq\qq+\(\int_r^Te^{-\b s}ds\)\(\int_r^Te^{\b s}\big|q(s,r,t)\big|^2ds\)\Big\}dtdr\\
\ns\ds={1\over\b}\dbE\int_\rho^T\int_\rho^r\int_\rho^s\(e^{\b(r+s)}-e^{\b(t+s)}\)|q(r,s,t)|^2dtdsdr\\
\ns\ds\qq+{1\over\b}\dbE\int_\rho^T\int_\rho^s\int_\rho^r\(e^{\b(t+s)}-e^{\b(r+s+t-T)}\)|q(s,r,t)|^2
dtdrds\\
\ns\ds\les{2\over\b}\dbE\int_\rho^T\int_\rho^r\int_\rho^s e^{\b(r+s)}\big|q(r,s,t)\big|^2dtdsdr\\
\ns\ds={2\over\b}\int_\rho^T\int_\rho^re^{\b(r+s)}\dbE\Big|\int_\rho^sq(r,s,t)dW(t)
\Big|^2dsdr\les{2\over\b}\dbE\int_\rho^T\int_\rho^re^{\b(r+s)}|p(r,s)|^2dsdr.\ea$$
To sum up the above arguments, with any $\rho\in[\t,T]$, we have
\bel{For-deriving-estimate}\ba{ll}
\ns\ds \dbE\int_{\rho}^T \int_{\rho}^r e^{\b(r+t)}| P(r,t)|^2 dt dr\\
\ns\ds \les  K \dbE\int_{\rho}^T \int_{\rho}^r e^{\b(r+t)}| F(r,t)|^2 dt dr +\frac{K}{\b}
 \dbE\int_{\rho}^T \int_{\rho}^r e^{\b (r+t)}|p(r,t)|^2  dt dr<\infty.
\ea
\ee
Hence
$P(\cd,\cd)\in L^2_\dbF(\Om;L^2(\D_*[\rho,T]))$, and one can define a map $\Xi$ from
$L^2_\dbF(\Om;L^2(\D_*[\rho,T]))$ to itself as $\Xi (p )=P$.

Suppose $\bar p $, $\wt p $ are two elements in $L^2_\dbF(\Om;L^2(\D_*[\rho,T]))$, and $(\bar q,\wt q)$ is defined similarly as the second equality in (\ref{Second-order-1-contraction}). By the previous arguments, we have $\bar P=\Xi(\bar p),$ $\wt P=\Xi(\wt p)$, and define $\bar Q$, $\wt Q$ accordingly. Also let
$$\h p=p-\wt p,\q\h q=q-\wt q,\q \h P=P-\wt P,\q\h Q=Q-\wt Q.$$
Then we have
\bel{Second-order-1-1-contraction}\left\{\ba{ll}
\ns\ds \h P(r,t)= \int_r^T\big[
b_x(s,t)^\top \dbE_t \h p(s,r)^{\top}+\si_x(s,t)^\top \h q(s,r,t)^{\top}\big]ds\\
\ns\ds\qq\qq+\int_t^r\big[b_x(s,t)^\top \dbE_t \h p(r,s)+\si_x(s,t)^\top \h q(r,s,t)\big]ds,\ \ r\ges t,\\
\ns\ds \dbE_\th \h p(r,t)=\h p(r,t)-\int_\th^{t} \h q(r,t,s)dW(s),\ \ \rho\les\th\les t\les r\les T.
\ea\right.
\ee
By the first equality in (\ref{Second-order-1-1-contraction}), for any $\rho\in[\t,T]$,
\bel{Lemma-useful-3}\ba{ll}
\ns\ds\dbE\int_\rho^T\int_\rho^re^{\b(r+t)}|\h P(r,t)|^2dtdr\\
\ns\ds\qq\les K\dbE\int_\rho^T\int_\rho^re^{\b(r+t)}\Big\{\(\int_t^r\big|\h p(r,s)\big|ds\)^2+\(\int_r^T\big|\h p(s,r)\big|ds\)^2\Big\}dtdr\\
\ns\ds\qq+K\dbE\int_\rho^T\int_\rho^re^{\b(r+t)}\Big\{\(\int_t^r\big|\h q(r,s,t)\big|ds\)^2
+\(\int_r^T\big|\h q(s,r,t)\big|ds\)^2\Big\}dtdr.\ea\ee
Similar as the estimation of \rf{estimate}, we have
for any $\rho\in[\t,T]$,
$$\dbE\int_\rho^T\int_\rho^re^{\b(r+t)}|\h P(r,t)|^2dtdr\les{K\over\b}\dbE\int_\rho^T \int_\rho^re^{\b(r+t)}|\h p(r,t)|^2dtdr.$$
Therefore, by choosing $\b>0$ large enough, we obtain the existence and uniqueness of $P$ on $[\rho,T]$, as well as the above conclusion (\ref{Estimate-1-preparation}).\endpf


\medskip

Now we return to system (\ref{Second-order adjoint equations}), the second-order adjoint system in our scenario.
\begin{theorem}\label{Well-posedness-SOAE}
\sl Suppose {\rm(H1)--(H2)} hold. Then \rf{Second-order adjoint equations}
admits a unique adapted M-solution $(P_i,Q_i)$, $1\les i\les 4$, such that $(P_3,P_4)\in L^2_{\dbF}(\t,T;\dbR^{n\times n})\times L^2_\dbF(\Om;L^2(\D_*[\t,T];\dbR^{n\times n}))$.
\end{theorem}

\it Proof. \rm By standard BSVIEs theory (\cite{Yong-2008}), the well-posedness in the spaces is easy to see
$$\ba{ll}
\ns\ds (P_1,Q_1)\in L^2_{\dbF}(\Omega;C([\t,T];\dbS^n))\times L^2_{\dbF}(\t,T;\dbS^n),\\
\ns\ds (P_2,Q_2)\in L^2_{\dbF}(\t,T;\dbR^{n\times n})\times L^2 (\t,T;L^2_{\dbF}(\t,T;\dbR^{n\times n})).\ea$$
We now prove the remaining well-posedness of $(P_3,P_4)$.
Recall (\ref{sE(e)*}), we see that $H_{xx}(\cd)\in L^2_{\dbF}(\t,T;\dbR^{n\times n})$.
By applying conditional expectation operator $\dbE_t$ on both sides of the equation of $(P_3,P_4)$ in \rf{Second-order adjoint equations}, we have
\bel{Third-fourth-equation-reduced}\left\{\2n\ba{ll}
\ns\ds P_3(t)=\BF_3(t)+\dbE_t\int_t^T\3n  \int_s^T\3n\(\si_x(\th,t)^\top P_4(\th,s)^\top\si_x(s,t)+\si_x(s,t)^\top\1n P_4(\th,s)\si_x(\th,t)\)d\th ds\\
\ns\ds\qq\qq+\dbE_t\int_t^T\3n\si_x(s,t)^\top P_3(s)\si_x(s,t)ds,\q
\t\les t\les T\\
\ns\ds P_4(r,t)\1n=\1n\BF_4(r,t)+\1n b_x(r,t)^\top\1n \dbE_t P_3(r)\1n+\1n\si_x(r,t)^\top Q_3(r,t)\\
\ns\ds\qq\qq\qq+\dbE_t\int_r^T\(b_x(s,t)^\top P_4(s,r)^{\top}+\si_x(s,t)^\top Q_4(s,r,t)^{\top}\)ds\\
\ns\ds\qq\qq\qq+\dbE_t\int_t^r\(b_x(s,t)^\top P_4(r,s)+\si_x(s,t)^\top Q_4(r,s,t)\)ds,\ \ \t\les t\les r\les T,\ea\right.\ee
where $\BF_3(\cd)$ and $\BF_4(\cd\,,\cd)$ are given by the following, depending on $(P_1,Q_1)$ and $(P_2,Q_2)$,
$$\ba{ll}
\ds\BF_3(t)\=H_{xx}(t)+\si_x(T,t)^\top P_1(t)\si_x(T,t)\\
\ds\qq\qq+\dbE_t\int_t^T\big[
 \si_x(s,t)^\top P_2(s)\si_x(T,t)+\si_x(T,t)^\top P_2(s)^\top\si_x(s,t)\big]ds,\\
\BF_4(r,t)\=b_x(T,t)^\top\dbE_tP_2(r)^\top+\si_x(T,t)^\top Q_2(r,t)^\top.\ea$$
Thus,
$$\dbE\(\int_\t^T|\BF_3(t)|^2dt+\int_\t^T\int_\t^r|\BF_4(r,t)|^2dtdr\)<\i.$$
For a given $p_3(\cd)\in L^2_{\dbF}(\t,T)\equiv L^2_{\dbF}(\t,T;\dbR^{n\times n})$, let the associated $q_3(\cd\,,\cd)$ be determined by martingale representation theorem:
$$p_3(t)=\dbE_rp_3(t)+\int_r^tq_3(t,s)dW(s),\qq\t\les r\les t\les T.$$
Then
\bel{p>q}\dbE|p_3(t)|^2=\dbE|\dbE_rp_3(t)|^2+\dbE\int_r^t|q_3(t,s)|^2ds\ges
\dbE\int_r^t|q_3(t,s)|^2ds.\ee
Now, let us consider
\bel{Step-2-main-theorem-1}\left\{\2n\ba{ll}
\ns\ds P_3(t)=\BF_3(t)\1n+\1n\dbE_t\int_t^T\3n\int_s^T\3n\(\si_x(\th,t)^\top P_4(\th,s)^\top\si_x(s,t)+\si_x(s,t)^\top\1n P_4(\th,s)\si_x(\th,t)\)d\th ds\\
\ns\ds\qq\qq+\dbE_t\int_t^T\3n\si_x(s,t)^\top p_3(s)\si_x(s,t)ds,\q
\t\les t\les T,\\
\ns\ds P_4(r,t)=\BF_4(r,t)+\1n b_x(r,t)^\top\1n \dbE_t p_3(r)\1n+\1n\si_x(r,t)^\top q_3(r,t)\\
\ns\ds\qq\qq+\dbE_t\int_r^T\(b_x(s,t)^\top P_4(s,r)^{\top}+\si_x(s,t)^\top Q_4(s,r,t)^{\top}\)ds\\
\ns\ds\qq\qq+\dbE_t\int_t^r\(b_x(s,t)^\top P_4(r,s)+\si_x(s,t)^\top Q_4(r,s,t)\)ds,\ \ \t\les t\les r\les T.\ea\right.\ee
By Lemma \ref{Preparation-1}, the second equation in \rf{Step-2-main-theorem-1} is solvable and
for any $\rho\in[\t,T]$ (noting \rf{p>q})
\bel{P_4<}\ba{ll}
\ns\ds\dbE\2n\int_\rho^T\3n\int_\rho^r\3n|P_4(r,t)|^2dtdr\1n\les\1n K\dbE\(\1n\int_\rho^T\3n\int_\rho^r\3n|\BF_4(r,t)|^2dtdr+\2n\int_\rho^T\3n|p_3(r)|^2dr
+\2n\int_\rho^T\3n\int_\rho^r \3n|q_3(r,t)|^2dtdr\)\\
\ns\ds\les K\dbE\(\int_\rho^T\int_\rho^r\3n|\BF_4(r,t)|^2dtdr+\int_\rho^T\3n |p_3(r)|^2dr\)<\i.\ea\ee
Therefore, for any given $\b>0$,
$$\ba{ll}
\ns\ds\dbE\int_\t^Te^{\b t}\Big|\int_t^T\int_s^T\si_x(s,t)^\top P_4(\th,s)\si_x(\th,t)d\th ds \Big|^2dt\\
\ns\ds\les K\dbE\int_\t^Te^{\b t}\(\int_t^T\int_s^T|P_4(\th,s)|^2d\th ds\)dt
=K\dbE\int_\t^Te^{\b t}\(\int_t^T\int_t^\th|P_4(\th,s)|^2dsd\th\)dt\\
\ns\ds\les K\dbE\int_\t^Te^{\b t}\(\int_t^T\int_t^r|\BF_4(r,t)|^2dtdr+\int_t^T|p_3(r)|^2dr\)dt\\
\ns\ds\les{K(e^{\b T}-e^{\b\t})\over\b}\dbE\int_\t^T\int_\t^r|\BF_4(r,t)|^2dtdr+K\dbE\int_\t^T\int_\t^re^{\b t}
|p_3(r)|^2dtdr\\
\ns\ds\les{Ke^{\b T}\over\b}\dbE\int_\t^T\int_\t^r|\BF_4(r,t)|^2dtdr+{K\over\b}\dbE\int_\t^Te^{\b r}
|p_3(r)|^2dtdr<\i.\ea$$
Similarly,
$$\ba{ll}
\ns\ds\dbE\int_\t^T e^{\b t}\Big|\int_t^T\int_s^T\si_x(\th,t)^\top P_4(\th,s)^\top
\si_x(s,t)d\th ds\Big|^2dt\\
\ns\ds\les{Ke^{\b T}\over\b}\dbE\int_\t^T\int_\t^r|\BF_4(r,t)|^2dtdr+{K\over\b}\dbE\int_\t^Te^{\b r}
|p_3(r)|^2dtdr<\i.\ea$$
Therefore, for the equality of $P_3(\cd)$ in (\ref{Step-2-main-theorem-1}), for any given $\b>0$, we have (noting \rf{p>q})
\bel{P_3<}\ba{ll}
\ns\ds\dbE\int_\t^T\3n e^{\b t}|P_3(t)|^2dt\les K\dbE\int_\t^T\3n e^{\b t}
\(|\BF_3(t)|^2+\int_t^T\3n\int_s^T\3n|P_4(\th,s)|^2d\th ds+\int_t^T|p_3(s)|^2ds\)dt\\
\ns\ds\les K\dbE\Big[e^{\b T}\int_\t^T|\BF_3(t)|^2dt+\int_\t^Te^{\b t}\(\int_t^T\int_t^\th|P_4(\th,s)|^2dsd\th+\int_t^T|p_3(s)|^2ds\)dt\Big]\\
\ns\ds\les K\dbE\Big[e^{\b T}\int_\t^T|\BF_3(t)|^2dt+{e^{\b T}\over\b}
\int_\t^T\int_\t^r|\BF_4(r,s)|^2dsdr+\int_\t^T\int_\t^se^{\b t}|p_3(s)|^2dtds\Big]\\
\ns\ds\les K\dbE\Big[e^{\b T}\int_\t^T|\BF_3(t)|^2dt+{e^{\b T}\over\b}
\int_\t^T\int_\t^r|\BF_4(r,s)|^2dsdr+{1\over\b}\int_\t^Te^{\b s}|p_3(s)|^2dtds\Big].\ea\ee
Hence, we have $P_3(\cd)\in L^2_{\dbF}(\t,T)$, and the map
$\Xi(p_3)=P_3$ in $L^2_{\dbF}(\t,T)$ is well-defined.

\ms

Now, suppose $\bar p_3$, $\wt p_3$ are two elements in $L^2_{\dbF}(\t,T)$. We then have
$\Xi(\bar p_3)=\bar P_3$, $\Xi(\wt p_3)=\wt P_3$, and define
$$\h p_3=\bar p_3-\wt p_3,\q\h P_3=\bar P_3-\wt P_3,\q\h P_4=\bar P_2-\wt P_4.$$
Then for any $t\in[\t,T]$, by \rf{P_4<}--\rf{P_3<} (with the corresponding $\BF_3(\cd)=0$ and $\BF_4(\cd\,,\cd)=0$),
$$\ba{ll}
\ns\ds \dbE\int_\t^T e^{\b r} |\h P_3(r)|^2 dr\les  \frac{K}{\b}\dbE
 \int_\t^T  e^{\b t}
 |\h p_3(t)|^2dt.\ea$$
By choosing $\b$ large, we see that the map $\Xi$ is a contraction. Hence, it admits a unique fixed point $P_3)\cd)$. The rest of the conclusion follows easily. \endpf

\begin{remark}
\rm
In the above two results, there are two crucial points, i.e., the introduction of equivalent $\b$-norm of $L^2_\dbF(\Om;L^2(\D_*[\t,T]))$ and $L^2_{\dbF}(\t,T)$, respectively, and the subtle use of Lemma \ref{Preparation-1} in Theorem \ref{Well-posedness-SOAE}. We point out that the idea of using multiplier $e^{\b t}$ ($\b$-argument, for short) appeared in the literature (e.g. \cite{Shi-Wang-2012-JKMS, Wang-Yong-2013, Shi-Wang-Yong-2015}).
\end{remark}
\medskip

\section{The Pontryagin's Type Maximum Principles and Related Extensions}

In this section, let us first give the main result of this paper, maximum principle of optimal controls $\bar u(\cd)$. Then we show that the result is consistent with the maximum principle of FSDEs. At last we present some extensions to multi-objective problem and multi-person differential games of governed by FSVIEs.

\subsection{Maximum principle for FSVIEs}

For convenience, let us recall the following first-order adjoint equations:
\bel{1st-adjoint-equation*}\left\{\2n\ba{ll}
\ns\ds\eta(t)=h_x^\top-\int_t^T\z(s)dW(s),\qq t\in[0,T],\\
\ns\ds Y(t)=g_x(t)^\top+b_x(T,t)^\top h_x^\top+\si_x(T,t)^{\top}\z(t)\\
\ns\ds\qq\q+\int_t^T\(b_x(s,t)^\top Y(s)+\si_x(s,t)^\top Z(s,t)\)ds-\int_t^TZ(t,s)dW(s),\q t\in[0,T],\ea\right.\ee
and the Hamiltonian
\bel{H*}\ba{ll}
\ns\ds H(s,x,u,\eta(s),\z(s),Y(\cd),Z(\cd\,,s))\1n=\1n\lan\eta(s),b(T,s,x,u)\ran
\1n+\1n\lan\z(s),\si(T,s,x,u)\ran\1n+g(s,x,u)\\
\ns\ds\qq\qq\qq\qq\qq\qq\qq\;+\dbE_s\2n\int_s^T\3n\(\2n\lan Y(t),b(t,s,x,u)\ran\1n+\1n\lan Z(t,s),\si(t,s,x,u)\ran\2n\)dt.\ea\ee
Then the Pontryagin maximum principle for Problem (C) can be stated as follows.

\begin{theorem}
\sl Let {\rm(H1)--(H2)} hold. Let $(\bar X(\cd).\bar u(\cd))$ be an optimal pair of Problem {\rm(C)}. Then
\bel{Maximum-principle}\ba{ll}
\ns\ds H(\t,\bar X(\t),u, \eta(\t),\z(\t), Y(\cd), Z(\cd\,,\t))
\!-\!H(\t,\bar X(\t),\bar u(\t), \eta(\t),\z(\t), Y(\cd), Z(\cd\,,\t))\\
\ns\ds \qq+ \frac 1 2 \Big\{\d\si(T,\t)^\top P_1(\t)\d\si(T,\t)+
\dbE_\t\int_\t^T\Big[\d\si(\th,\t)^\top P_3(\th)\d\si(\th,\t)\\
\ns\ds\qq\qq+\d\si(T,\t)^\top P_2(\th)^\top\d\si(\th,\t)+\d\si(\th,\t)^\top
P_2(\th)\d\si(T,\t)\\
\ns\ds\qq+ \int_\th^T\big[\d\si(\th',\t)^\top P_4(\th',\th)^{\top}
 \d\si
(\th,\t)+\d\si(\th,\t)^\top P_4(\th',\th)
 \d\si
(\th',\t)\big]d\th'\Big]d\th\Big\}\ges0,\\
\ns\ds\qq\qq\qq\qq\qq\qq\qq\qq\qq\qq\as\2n,~\ae\t\in[0,T],~\forall u\in U,\ea\ee
where $H$ is defined in \rf{H*}, $(\eta(\cd),\z(\cd),\bar Y(\cd),\bar Z(\cd))$ satisfies BSVIE \rf{1st-adjoint-equation*}, and $(P_i(\cd),Q_i(\cd))$ ($i=1,2,3,4$) is the unique adapted M-solution of the second order-adjoint equation \rf{Second-order adjoint equations}.
\end{theorem}

\it Proof. \rm Recall that
\bel{Delta-sigma-small}\ba{ll}
\ns\ds \d\si^\e(t,s)\=\big[\si(t,s,\bar X(s), u)-\si(t,s,\bar X(s), \bar u(s))\big]{\bf1}_{[\t,\t+\e]}(s).
\ea
\ee
Thanks to  (\ref{Delta-sigma-small}) and Lebesgue differentiability theorem, for almost $\t\in[0,T]$, by Theorem \ref{variation-2}
\bel{sE(e)*-final-maximum}\ba{ll}
\ns\ds \lim_{\e\to 0}\frac{\sE(\e)}{\e}=\dbE\Big\{\d\si(T,\t)^\top P_1(\t)\d\si(T,\t)+
\dbE_\t\int_\t^T\Big[\d\si(T,\t)^\top
P_2(\th)^\top\d\si(\th,\t)\\
\ns\ds\qq+\d\si(\th,\t)^\top
P_2(\th)\d\si(T,\t)+
\d\si(\th,\t)^\top P_3(\th)\d\si(\th,\t)d\th
\\
\ns\ds\qq+ \int_\th^T\big[\d\si(\th',\t)^\top P_4(\th',\th)^{\top}
 \d\si
(\th,\t)+\d\si(\th,\t)^\top P_4(\th',\th)
 \d\si
(\th',\t)\big]d\th'\Big]d\th\Big\}.\ea\ee
According to Theorems \ref{variation} and \ref{Well-posedness-SOAE}, by the optimality of $\bar u(\cd)$ and the arbitrariness of $u\in U$, we obtain the above maximum condition (\ref{Maximum-principle}) immediately.\endpf

\subsection{The case of FSDEs}

In this subsection, we show that the above maximum principle reduces to the FSDEs case presented in \cite{Peng-1990, Yong-Zhou-1999}, when the following holds:
\bel{b-si}b(t,s,x,u)=b(s,x,u),\qq\si(t,s,x,u)=\si(s,x,u),\qq\f(t)=x.\ee

Let us look at the second-order adjoint system (\ref{Second-order adjoint equations}) under \rf{b-si}. The first equation is unchanged. The second equation becomes
$$\ba{ll}
\ns\ds P_2(t)\1n=\1n b_x(t)^{\1n\top}\1n P_1(t)\1n+\1n\si_x(t)^{\1n\top}\1n Q_1(t)\1n+\2n\int_t^T\3n\(b_x(t)^{\1n\top}\1n P_2(s)\1n+\1n\si_x(t)^{\1n\top}\1n Q_2(s,t)\)ds-\int_t^TQ_2(t,s)dW(s)\\
\ns\ds\qq=b_x(t)^\top\(P_1(t)+\int_t^TP_2(s)ds\)+\si_x(t)^\top\(Q_1(t)+\int_t^T
Q_2(s,t)ds\)-\int_t^TQ_2(t,s)dW(s)\\
\ns\ds\equiv b_x(t)^\top\cP_2(t)+\si_x(t)^\top\cQ_2(t)-\int_t^TQ_2(t,s)dW(s).\ea$$
The third equation becomes
$$\ba{ll}
\ns\ds P_3(t)=H_{xx}(t)+\si_x(t)^\top P_1(t)\si_x(t)+\int_t^T\3n\Big[\si_x(t)^\top P_2(s)^\top\si_x(t)+\si_x(t)^\top
P_2(s)\si_x(t)\\
\ns\ds\qq\qq+\int_t^T\3n\(\si_x(t)^\top P_4(\th,s)^\top\si_x(t)+\si_x(t)^\top\1n P_4(\th,s)\si_x(\t)\)d\th\Big]ds\\
\ns\ds\qq\qq+\int_t^T\3n\si_x(t)^\top P_3(s)\si_x(t)ds-\2n\int_t^TQ_3(t,s)dW(s)\\
\ns\ds=H_{xx}(t)+\si_x(t)^\top\Big\{P_1(t)+\int_t^T\3n\Big[P_2(s)^\top+P_2(s)+P_3(s)+\int_s^T\3n
\(P_4(\th,s)^\top+P_4(\th,s))\)d\th\Big]ds\Big\}\si_x(t)\\
\ns\ds\qq\qq\qq\qq-\int_t^TQ_3(t,s)dW(s)\\
\ns\ds\equiv H_{xx}(t)+\si_x(t)^\top\cP_3(t)\si_x(t)-\int_t^TQ_3(t,s)dW(s).\ea$$
The fourth equation becomes
$$\ba{ll}
\ns\ds P_4(r,t)\1n=\1n b_x(t)^\top\1n P_2(r)^\top\3n+\1n\si_x(t)^\top\1n Q_2(r,t)^\top\3n+\1n b_x(t)^\top\1n P_3(r)\1n+\1n\si_x(t)^\top Q_3(r,t)\\
\ns\ds\qq\qq\qq+\int_r^T\(b_x(t)^\top P_4(s,r)^\top+\si_x(t)^\top Q_4(s,r,t)^\top\)ds\\
\ns\ds\qq\qq\qq+\int_t^r\(b_x(t)^\top P_4(r,s)+\si_x(t)^\top Q_4(r,s,t)\)ds-\int_t^TQ_4(r,t,s)dW(s)\\
\ns\ds=\1n b_x(t)^\top\(P_2(r)^\top+P_3(r)+\int_r^TP_4(s,r)^\top ds+\int_t^rP_4(r,s)ds\)\\
\ns\ds\q+\si_x(t)^\top\(Q_2(r,t)^\top+Q_3(r,t)+\int_r^TQ_4(s,r,t)^\top ds+\int_t^rQ_4(r,s,t)ds\)-\int_t^TQ_4(r,t,s)dW(s)\\
\ns\ds\equiv b_x(t)^\top\cP_4(r,t)+\si_x(t)^\top\cQ_4(r,t)-\int_t^TQ_4(r,t,s)dW(s),\ea$$
where
$$\left\{\2n\ba{ll}
\ns\ds\cP_2(t)=P_1(t)+\int_t^TP_2(s)ds,\qq\cQ_2(t)=Q_1(t)+\int_t^TQ_2(s,t)ds,\\
\ns\ds\cP_3(t)=P_1(t)+\int_t^T\3n\Big[P_2(s)^\top+P_2(s)+P_3(s)+\int_s^T\3n
\(P_4(\th,s)^\top+P_4(\th,s))\)d\th\Big]ds,\\
\ns\ds\cP_4(r,t)=P_2(r)^\top+P_3(r)+\int_r^TP_4(s,r)^\top ds+\int_t^rP_4(r,s)ds,\\
\ns\ds\cQ_4(r,t)=Q_2(r,t)^\top+Q_3(r,t)+\int_r^TQ_4(s,r,t)^\top ds+\int_t^rQ_4(r,s,t)ds.\ea\right.$$
The Hamiltonian takes the form
\bel{H*}\ba{ll}
\ns\ds H(s,x,u,\eta(s),\z(s),Y(\cd),Z(\cd\,,s))\\
\ns\ds=\1n\lan\eta(s)+\dbE_s\int_s^TY(t)dt,b(s,x,u)\ran
\1n+\1n\lan\z(s)+\dbE_s\int_s^TZ(t,s)dt,\si(s,x,u)\ran\1n+g(s,x,u)\\
\ns\ds\equiv\lan\cY(s),b(s,x,u)\ran+\lan\cZ(s),\si(s,x.u)\ran+g(s,x,u)\equiv
\cH(s,x,u,\cY(s),\cZ(s)),\ea\ee
where
\bel{cY}\cY(s)=\eta(s)+\dbE_s\int_s^TY(t)dt,\qq\cZ(s)=\z(s)+\dbE_s\int_s^T
Z(t,s)dt,\ee
and
\bel{cH}\cH(s,x,u,\cY,\cZ)=\lan\cY,b(s,x,u)\ran+\lan\cZ,\si(s,x,u)\ran+g(s,x,u).\ee
The maximum condition becomes
\bel{Maximum-principle*}\ba{ll}
\ns\ds0\les\cH(\t,\bar X(\t),u,\cY(\t),\cZ(\t))
\!-\!H(\t,\bar X(\t),\bar u(\t),\cY(\t),\cZ(\t))\\
\ns\ds \qq+{1\over2}\d\si(\t)^\top\Big[P_1(\t)+
\dbE_\t\int_\t^T\(P_2(\th)^\top+P_2(\th)+P_3(\th)\\
\ns\ds\qq\qq\qq\qq+\int_\th^T\big[ P_4(\th',\th)^\top
+P_4(\th',\th)\big]d\th'\)d\th\Big]\d\si(\t)\\
\ns\ds=\cH(\t,\bar X(\t),u,\cY(\t),\cZ(\t))
\!-\!H(\t,\bar X(\t),\bar u(\t),\cY(\t),\cZ(\t))+{1\over2}\d\si(\t)^\top\dbE_\t\cP_3(\t)\d\si(\t)\ea\ee
Next we check that $(\cY,\cZ)$ satisfies the following first-order adjoint equation
\bel{first-order-adjoint-equation-SDEs-0}\left\{\ba{ll}
\ns\ds d\cY(t)=-\Big[g_x(t)^{\top}+  b_x(t)^{\top}\cY(t)+ \si_x(t)^{\top}\cZ(t)\Big]dt+\cZ(t)dW(t),\ \ t\in[0,T],\\
\ns\ds \cY(T)=h_x(\bar X(T))^{\top}.
\ea\right.\ee
Note that in the current case, we have
$$\ba{ll}
\ns\ds Y(t)=g_x(t)^\top+b_x(t)^\top h_x^\top+\si_x(t)^{\top}\z(t)\\
\ns\ds\qq\qq+\int_t^T\(b_x(t)^\top Y(s)+\si_x(t)^\top Z(s,t)\)ds-\int_t^TZ(t,s)dW(s)\\
\ns\ds\qq=g_x(t)^\top+b_x(t)^\top\(h_x^\top+\int_t^TY(s)ds\)
+\si_x(t)^\top\(\z(t)+\int_t^TZ(s,t)ds\)\\
\ns\ds\qq\qq -\int_t^TZ(t,s)dW(s).
\ea$$
Thus, applying $\dbE_t$ on both sides, we get
\bel{First-order-adjoint-proof-1}\ba{ll}
\ns\ds Y(t)=g_x(t)^\top+b_x(t)^\top\(\eta(t)+\dbE_t\int_t^TY(s)ds\)+\si_x(t)^\top\(\z(t)
+\int_t^TZ(s,t)ds\)\\
\ns\ds \qq=g_x(t)^\top+b_x(t)^\top \cY(t)+\si_x(t)^\top\cZ(t),
\ea\ee
By Fubini theorem,
$$\ba{ll}
\ns\ds \int_t^T Y(s)ds=\dbE_t\int_t^T Y(s)ds+\int_t^T \int_t^s Z(s,r)dW(r)ds\\
\ns\ds\qq\qq\q=\dbE_t\int_t^T Y(s)ds+\int_t^T \int_r^T Z(s,r) ds dW(r).\ea$$
Hence
$$\ba{ll}
\ns\ds h_x^{\top}+ \int_t^T Y(s)ds=\eta(t)+\dbE_t\int_t^T Y(s)ds+\int_t^T
\Big[\int_r^T Z(s,r) ds+\zeta(r)\Big]dW(r)\\
\ns\ds \qq\qq\qq\qq =\cY(t)+\int_t^T \cZ(r)dW(r).
\ea
$$
Plugging (\ref{First-order-adjoint-proof-1}) into the above equation, we end up with BSDE (\ref{first-order-adjoint-equation-SDEs-0}), the first-order adjoint equation.

\ms

Finally, it suffices to prove that $\sM(\cd)\=\dbE_{\cd} \cP_3(\cd)$ satisfies the second order-adjoint equation
\bel{Second-order-SDEs-case}\ba{ll}
\ns\ds \sM(r)=h_{xx}(\bar X(T))+\int_r^T \Big[ b_x(t)^{\top}\sM(t)+ \si_x(t)^{\top}\sN(t)+\sM(t) b_x(t) + \sN(t) \si_x(t)\\
\ns\ds\qq\qq  +  H_{xx}(t)
+ \si_x(t)^{\top}\sM(t)  \si_x(t)\Big]dt-\int_r^T\sN(t) dW(t),\ \ r\in[0,T],
\ea
\ee
where
\bel{Notations-sM-sN-1}\ba{ll}
\ns\ds \sN(t)\=\cQ_2(t) +\int_t^T \cQ_4(s,t)ds.\ea\ee
To this end, we observe the following: For the equation of $P_2$, one has
$$\ba{ll}
\ns\ds\int_r^T\big[P_2(s)+P_2(s)^\top\big]ds\\
\ns\ds=\int_r^T\Big\{b_x(s)^\top\dbE_s\cP_2(s)+\si_x(s)^{\top}\cQ_2(s)+\dbE_s\cP_2(s)^\top b_x(s)+\cQ_2(s)^\top\si_x(s)\Big\}ds.\ea$$
For the equation of $P_3$, one has
%
$$\int_r^T P_3(t)dt  =
\int_r^T   H_{xx}(t)dt + \int_r^T  \si_x(t)^{\top}\dbE_t\cP_3(t)\si_x(t)dt.$$
For the equation of $P_4$, we have
%
%
%
$$\int_r^T\int_t^T P_4(s,t)dsdt =
\int_r^T\int_t^T \Big[  b_x(t)^{\top}\dbE_s\cP_4(s,t) +
\si_x(t)^{\top}\cQ_4(s,t)\Big]dsdt.$$
Therefore,
\bel{sG-equality-expression}\ba{ll}
\ns\ds \cP_3(r)
=h_{xx}(\bar X(T))+\int_r^T \Big[ b_x(t)^{\top}\sM(t)+ \si_x(t)^{\top}\sN(t)+\sM(t)  b_x(t)\\
\ns\ds\qq \qq\q + \sN(t) \si_x(t)+  H_{xx}(t)
+ \si_x(t)^{\top}\sM(t)\si_x(t)\Big]dt,
\ea
\ee
where $\sN(\cd)$ is defined in (\ref{Notations-sM-sN-1}).

By the definition of $Q_2$, $Q_3$, and Fubini theorem, we see that
\bel{Martingale-Part-III}\ba{ll}
\ns\ds \int_r^T P_i(t)dt =\dbE_r \int_r^T P_i(t)dt+\int_r^T\int_s^T Q_i(t,s)dt dW(s), \ \  i=2,3.
\ea
\ee
Using Fubini theorem, we have
$$\ba{ll}
\ns\ds  \int_r^T \int_t^T \int_r^t
Q_4(s,t,\th)dW(\th)dsdt=\int_r^T \int_\th^T \int_t^T
Q_4(s,t,\th)ds dt dW(\th).
\ea
$$
Therefore,
\bel{Martingale-Part-I}\ba{ll}
\ns\ds \int_r^T\3n\int_t^T \big[P_4(s,t)+P_4(s,t)^{\top}\big] dsdt\\
\ns\ds =\dbE_r\2n\int_r^T\3n\int_t^T\3n\big[P_4(s,t)\1n+\1n P_4(s,t)^{\top}\big]dsdt\1n+\2n\int_r^T\3n \int_\th^T\3n\int_t^T\3n\big[
Q_4(s,t,\th)+Q_4(s,t,\th)^{\top}\big]ds dt dW(\th).\ea\ee
Combining (\ref{Martingale-Part-III}) and (\ref{Martingale-Part-I}), we have
$$\ba{ll}
\ns\ds \cP_3(r)  =\sM(r)+ \int_r^T\sN(\th) dW(\th).\ea$$
Consequently, by means of (\ref{sG-equality-expression}), for $r\in[0,T],$ we have BSDE (\ref{Second-order-SDEs-case}), the second-order adjoint equation.

To sum up, the maximum condition \rf{Maximum-principle*}, the first-order adjoint equation \rf{first-order-adjoint-equation-SDEs-0} and the second-order adjoint equation \rf{Second-order-SDEs-case} form the Pontryagin's maximum principle for optimal control of  FSDEs.

\subsection{Multi-person dynamic games for FSVIEs}

For integer $N\ges2$, we consider an $N$-person dynamic game for an FSVIE. In this case, the cost functional is defined as in (\ref{cost3.1}), where $h:\Om\times\dbR^n\rightarrow\dbR^N$, $g:[0,T]\times\Om\times\dbR^n\times\dbR^m\to\dbR^N$ are vector-valued functions.
For simplicity, we only discuss the non-cooperative dynamic games, namely, the $\ell$-th player in the game wants to minimize his/her own cost functional $J^\ell(u(\cd))$, regardless of other players' cost functional.

Let $u(\cd)=(u_1(\cd),u_2(\cd),\cdots, u_N(\cd))$ with $u_\ell(\cd)\in\cU^\ell_{ad}$. Here $\cU_{ad}^\ell$ is defined associated with $U_\ell\in\dbR^{m_\ell}$.
For notational simplicity, let
$$(u^c_\ell,v)\=(u_1,\cds,u_{\ell-1},v,u_{\ell+1},\cds,u_N),\q 1\les\ell\les N.$$
Then player $\ell$ selects $u_\ell(\cd)\in\cU_{ad}^\ell$ to minimize the functional
$$v(\cd)\mapsto J^\ell(u^c_\ell(\cd),v(\cd))\equiv J^\ell(u_1(\cd),\cds,u_{\ell-1},v(\cd),u_{\ell+1},\cds,u_N(\cd)).$$
Obviously, $J^\ell(u(\cd))$ not only depends on $u_\ell(\cd)$, but also $u_k(\cd)$, $k\ne\ell$. Therefore, the optimal control of Player $\ell$ depends on the controls of the other players.

\begin{definition} \rm An $N$-tuple $\bar u(\cd)=(\bar u_1(\cd),\cdots,\bar u_N(\cd))\in \prod_{\ell=1}^N \cU_{ad}^\ell$ is called an open-loop Nash equilibrium of the game if the following holds:
$$J^\ell(\bar u(\cd))\les J^\ell(\bar u^c_\ell(\cd),v(\cd)),\qq\forall v(\cd)\in\cU_{ad}^\ell,~
1\les\ell\les N.$$
\end{definition}

We have the following Pontryagin type maximum principle for Nash equilibria of the $N$-person dynamic game of FSVIEs. To state the result, let $(\eta^\ell,\z^\ell,Y^\ell,Z^\ell)$ be the solution of (\ref{1st-adjoint-equation}) associated with $(h_x^\ell(\bar X(T)), g_x^\ell(\cd))$. Then we define $H^\ell(\cd)$ as in (\ref{H}) accordingly. Similarly, with $(h_{xx}^\ell, H_{xx}^\ell)$, let $(P_1^\ell,P_2^\ell,P_3^\ell,P_4^\ell)$ be the solution of (\ref{Second-order adjoint equations}).

\begin{theorem}
\sl Let {\rm(H1)--(H2)} hold with $g$ and $h$ being $\dbR^N$-valued. Let $\bar u=(\bar u_1(\cd),\cdots,\bar u_N(\cd))\in\prod_{\ell=1}^N\cU_{ad}^\ell$ be a Nash equilibrium of the game. Then for any $u\in U^\ell$,
\bel{Maximum-principle-Game-1}\ba{ll}
\ns\ds H^\ell(\t,\bar X(\t),\bar u^c_\ell(\t),u,\eta^\ell(\t),\z^\ell(\t), Y^\ell(\cd), Z^\ell(\cd\,,\t))\\
\ns\ds\qq -\!H^\ell(\t,\bar X(\t),\bar u(\t),\eta^\ell(\t),\z^\ell(\t),Y^\ell(\cd), Z^\ell(\cd\,,\t))\\
\ns\ds\qq+{1\over2}\Big\{\d\si(T,\t)^\top P_1^\ell(\t)\d\si(T,\t)+
\dbE_\t\int_\t^T\Big[\d\si(T,\t)^\top
P_2^\ell(\th)^\top\d\si(\th,\t)\\
\ns\ds\qq+\d\si(\th,\t)^\top P_2^\ell(\th)\d\si(T,\t)+
\d\si(\th,\t)^\top P_3^\ell(\th)\d\si(\th,\t)d\th\\
\ns\ds\qq+ \int_\th^T\big[\d\si(\th',\t)^\top P_4^\ell(\th',\th)^\top\d\si
(\th,\t)+\d\si(\th,\t)^\top P_4^\ell(\th',\th)\d\si
(\th',\t)\big]d\th'\Big]d\th\Big\} \ges 0.\ea\ee
\end{theorem}

\ms

Let us continue to look at a special case: the two-person zero-sum dynamical games.
In this case, $N=2$, and $J^1(\bar u(\cd))+J^2(\bar u(\cd))=0$. Now, if $(\bar u_1(\cd),\bar u_2(\cd))\in \cU_{ad}^1\times \cU_{ad}^2$ is a Nash equilibrium, then
$$\ba{ll}
\ns\ds J^1(\bar u_1(\cd),\bar u_2(\cd))\les J^1(u_1(\cd),\bar u_2(\cd)),\qq\forall u_1(\cd)\in\cU_{ad}^1,\\
\ns\ds J^2(\bar u_1(\cd),\bar u_2(\cd))\les J^2(\bar u_1(\cd), u_2(\cd)),\qq\forall u_2(\cd)\in\cU_{ad}^2.\ea$$
This also implies that (denoting $\bar J(u(\cd))=J^1(u(\cd))$)
$$\bar J(\bar u_1(\cd),u_2(\cd))\les  \bar J (\bar u_1(\cd),\bar u_2(\cd))\les  \bar J (  u_1(\cd),\bar u_2(\cd)).$$
Hence $(\bar u_1(\cd),\bar u_2(\cd))$ is referred to as a {\it saddle point} of the game. In this case, to state the maximum principle, let $(\eta,\zeta,Y,Z,P_i)$ satisfies (\ref{1st-adjoint-equation}), (\ref{Second-order adjoint equations}) associated with $(h^1,g^1)$. We define
\bel{Maximum-principle-Game-2}\ba{ll}
\ns\ds\sH(\t,u_1,u_2)\=H(\t,\bar X(\t),u_1, u_2, \eta(\t),\z(\t), Y(\cd), Z(\cd,\t))\\
\ns\ds\qq={1\over2}\Big\{\d\si(T,\t)^\top P_1(\t)\d\si(T,\t)+\dbE_\t\int_\t^T\Big[\d\si(T,\t)^\top
P_2(\th)^\top\d\si(\th,\t)\\
\ns\ds\qq\q+\d\si(\th,\t)^\top P_2(\th)\d\si(T,\t)+\d\si(\th,\t)^\top P_3(\th)\d\si(\th,\t)d\th\\
\ns\ds\q\qq+ \int_\th^T\big[\d\si(\th',\t)^\top P_4(\th',\th)^{\top}
 \d\si(\th,\t)+\d\si(\th,\t)^\top P_4(\th',\th)\d\si(\th',\t)\big]d\th'\Big]d\th\Big\}.\ea\ee
\begin{theorem}
\sl Let {\rm(H1)--(H2)} hold with $g$ and $h$ being $\dbR^2$-valued. Suppose $(\bar u_1(\cd),\bar u_2(\cd))\in\cU_{ad}^1\times\cU_{ad}^2$ is the a saddle point of the game. Then for any $(u_1,u_2)\in U^1\times U^2$,
$$\sH(\t,\bar u_1(\t),u_2)\les \sH(\t,\bar u_1(\t),\bar u_2)\les \sH(\t, u_1(\t),\bar u_2).$$
\end{theorem}

%
%
%

\section{A Proof of Theorem 3.1.}

We would like to carry out the variations of the state process and the cost functional under the spike variation, The main idea is the same as that for SDEs (see \cite{Peng-1990,Yong-Zhou-1999}). For convenience, let us recall the following abbreviations:
\bel{abb}\left\{\2n\ba{ll}
\ns\ds b_x(t,s)=b_x(t,s,\bar X(s),\bar u(s)),\qq b_{xx}^i(t,s)=b_{xx}^i(t,s,\bar X(s),\bar u(s)),\qq1\les i\les n,\\
\ns\ds\si_x(t,s)=\si_x(t,s,\bar X(s),\bar u(s)),\qq\si_{xx}^i(t,s)=\si_{xx}^i(t,s,\bar X(s),\bar u(s)),\qq1\les i\les n,\\
\ns\ds g_x(s)=g_x(s,\bar X(s),\bar u(s)),\qq g_{xx}(s)=g_{xx}(s,\bar X(s),\bar u(s)),\\
\ns\ds h_x=h_x(\bar X(T)),\qq h_{xx}=h_{xx}(\bar X(T)),\\
\ns\ds\d b(t,s)=b(t,s,\bar X(s),u)-b(t,s,\bar X(s),\bar u(s)),\\
\ns\ds\d\si(t,s)=\si(t,s,\bar X(s),u)-\si(t,s,\bar X(s),\bar u(s)),\\
\ns\ds\d\si_x(t,s)=\si_x(t,s,\bar X(s),u)-\si_x(t,s,\bar X(s),\bar u(s))\ea\right.\ee
and the following variational equations:
\bel{X_1-equation*} X_1^\e(t)\1n=\2n\int_0^t\3n b_x(t,s)X_1^\e(s)ds+\2n\int_0^t\2n\big[\si_x(t,s)X_1^\e(s)+\d\si(t,s){\bf1}_{[\t,\t+\e]}(s)
\big]dW(s),\q t\in[0,T],\ee
\bel{X_2-equation*}\ba{ll}
\ss\ds X_2^\e(t)=\int_0^t\(b_x(t,s)X_2^\e(s)+{1\over2}b_{xx}(t,s)X_1^\e(s)^2+\d b(t,s){\bf1}_{[\t,\t+\e]}(s)\)ds\\
\ss\ds\qq\qq+\int_0^t\(\si_x(t,s)X_2^\e(s)+{1\over2}\si_{xx}(t,s)X_1^\e(s)^2+\d\si_x(t,s)X_1^\e(s){\bf1}_{[\t,\t+\e]}(s)
\)dW(s),\\
\ns\ds\qq\qq\qq\qq\qq\qq\qq\qq\qq\qq\qq\qq\qq\qq t\in[0,T],\ea\ee
Our calculation is split into several steps.

\ms

\it Step 1. Estimation of $X^\e(\cd)-\bar X(\cd)$. \rm

\ms

Let us first look at the following
\bel{X-X1}\ba{ll}
\ds X^\e(t)-\bar X(t)=\int_0^t\(b\big(t,s,X^\e(s),u^\e(s)\big)-b\big(t,s,\bar X(s),\bar u(s)\big)\)ds\\
\ds\qq\qq\qq\qq+\int_0^t\(\si\big(t,s,X^\e(s),u^\e(s)\big)-\si\big(t,s,
\bar X(s),\bar u(s)\big)\)dW(s)\\
\ds\qq\qq\qq=\2n\int_0^t\(\wt b_x(t,s)[X^\e(s)\1n-\1n\bar X(s)]\1n+\1n\d b(t,s){\bf1}_{[\t,\t\1n+\e]}(s)\)ds\\
\ds\qq\qq\qq\qq+\int_0^t\(\wt\si_x(t,s)[X^\e(s)\1n-\1n\bar X(s)]\1n+\1n\d \si(t,s){\bf1}_{[\t,\t\1n+\e]}(s)\)dW(s),\ea\ee
where $\d b(t,s)$ and $\d\si(t,s)$ are given by \rf{abb}, and
$$\left\{\2n\ba{ll}
\ns\ds\wt b_x(t,s)=\int_0^1b_x\big(t,s,\bar X(s)\1n+\1n\a[X^\e(s)\1n-\1n\bar X(s)],u^\e(s)\big)d\a,\\
\ns\ds\wt\si_x(t,s)=\int_0^1\si_x\big(t,s,\bar X(s)\1n+\1n\a[X^\e(s)\1n-\1n\bar X(s)],u^\e(s)\big)d\a,\ea\right.$$
By the boundedness of $\wt b_x(\cd\,,\cd)$ and $\wt\si_x(\cd\,,\cd)$ (see (H1)), as well as the stability estimate of FSVIE (see Proposition 2.1), for $1\les k\les p$, we have
$$\ba{ll}
\ns\ds\sup_{t\in[0,T]}\dbE|X^\e(t)-\bar X(t)|^k\les K\dbE\Big[\(\int_\t^{\t+\e}
|\d b(t,s)|ds\)^k+\(\int_\t^{\t+\e}|\d\si(t,s)|^2ds\)^{k\over2}\Big]\\ [3mm]
\ss\ds\qq\qq\qq\qq\qq\les K\Big[\dbE\(\int_\t^{\t+\e}\3n\big(1+|\bar X(s)|\big)ds\)^k+\dbE\(\int_\t^{\t+\e}\3n\big(1+|\bar X(s)|^2\big)ds\)^{k\over2}\Big]\\
\ss\ds\qq\qq\qq\qq\qq\les K\(\e^{k-1}\int_\t^{\t+\e}(1+\dbE|\bar X(s)|^k)ds+\e^{k-2\over2}\int_\t^{\t+\e}(1+\dbE|\bar X(s)|^k)ds\)\\
\ss\ds\qq\qq\qq\qq\qq\les K \big[\e^k+\e^{k\over2}\big]\(1+\sup_{s\in[0,T]}\dbE|\bar X(s)|^k\)\les K\e^{k\over2}.\ea$$
This gives the first relation in \rf{|X_1|,|X_2|}.

\ms

\it Step 2. Estimation of $X^\e(\cd)-\bar X(\cd)-X_1^\e(\cd)$. \rm

\ms

We observe that
$$\ba{ll}
\ss\ds b(t,s,X^\e(s),u^\e(s))-b(t,s,\bar X(s),\bar u(s))\\
\ss\ds=\(\int_0^1b_x\big(t,s,\bar X(s)+\a[X^\e(s)-\bar X(s)],u^\e(s)\big)d\a\)
[X^\e(s)-\bar X(s)]+\d b(t,s){\bf1}_{[\t,\t+\e]}(s)\\
\ns\ds=b_x(t,s)[X^\e(s)-\bar X(s)]+\(\d b(t,s)+\D_\e^1 b(t,s)[X^\e(s)-\bar X(s)]\){\bf1}_{[\t,\t+\e]}(s),\ea$$
where $b_x(t,s)$ and $\d b(t,s)$ are as in \rf{abb}, and
$$\ds\D_\e^1b(t,s)\1n=\2n\int_0^1\3n\(b_x\big(t,s,\bar X(s)+\a[X^\e(s)-\bar X(s)],u^\e(s)\big)-b_x(t,s)\)d\a,\q t\in[0,T].$$
Likewise,
$$\ba{ll}
\ss\ds\si(t,s,X^\e(s),u^\e(s))-\si(t,s,\bar X(s),\bar u(s))\\
\ns\ds=\si_x(t,s)[X^\e(s)-\bar X(s)]+\(\d\si(t,s)+\D_\e^1\si(t,s)[X^\e(s)-\bar X(s)]\){\bf1}_{\t,\t+\e]}(s),\ea$$
where $\si_x(t,s)$ and $\d\si(t,s)$ are as in \rf{abb}, and
$$\D_\e^1\si(t,s)\1n=\2n\int_0^1\1n\(\si_x\big(t,s,\bar X(s)+\a[X^\e(s)-\bar X(s)],u^\e(s)\big)-\si_x(t,s)\)d\a,\q t\in[0,T].$$
%
%
%
Thus, one may also rewrite \rf{X-X1} as follows
\bel{X-X2}\ba{ll}
\ds X^\e(t)-\bar X(t)\\
\ds=\int_0^t\Big[b_x(t,s)[X^\e(s)-\bar X(s)]+\(\d b(t,s)+\D_\e^1 b(t,s)[X^\e(s)-\bar X(s)]\){\bf1}_{[\t,\t+\e]}(s)\Big]ds\\
\ds\q+\2n\int_0^t\2n\Big[\si_x(t,s)[X^\e(s)\1n-\1n\bar X(s)]\1n+\1n\(\d\si(t,s)\1n+\1n\D_\e^1\si(t,s)[X^\e(s)\1n-\1n\bar X(s)]\){\bf1}_{[\t,\t+\e]}(s)\Big]dW(s).\ea\ee
Consequently, we have
\bel{X-X-X1}\ba{ll}
\ds X^\e(t)-\bar X(t)-X^\e_1(t)\\
\ds=\2n\int_0^t\Big[b_x(t,s)[X^\e(s)\1n-\1n\bar X(s)\1n-\1n X_1^\e(s)]\1n+\1n\big(\d b(t,s)\1n+\1n\D_\e^1 b(t,s)[X^\e(s)\1n-\1n\bar X(s)]\big){\bf1}_{[\t,\t+\e]}(s)\Big]ds\\
\ds\q+\2n\int_0^t\2n\Big[\si_x(t,s)[X^\e(s)\1n-\1n\bar X(s)\1n-\1n X_1^\e(s)]\1n+\1n\big(\D_\e^1\si(t,s)[X^\e(s)\1n-\1n\bar X(s)]\big){\bf1}_{[\t,\t+\e]}(s)\Big]dW(s).\ea\ee
Note that $X_1^\e(\cd)$, the solution of \rf{X_1-equation*} satisfies (for any $1\les k\les p$)
\bel{|X_1}\ba{ll}
\ds\sup_{t\in[0,T]}\dbE|X^\e_1(t)|^k\les K\sup_{t\in[0,T]}\dbE\(\int_0^t\3n|\d\si(t,s)|^2{\bf1}_{[\t,\t+\e]}(s)ds
\)^{k\over2}\\
\ds\qq\qq\qq\q\les K\dbE\(\int_\t^{\t+\e}\big(1+|\bar X(s)|^2\big)ds\)^{k\over2}\les K\e^{k\over2},\ea\ee
which is the second relation \rf{|X_1|,|X_2|}. Hence, for $1\les k\les p$,
\bel{|X-X_1|}\ba{ll}
\ds\sup_{t\in[0,T]}\dbE|X^\e(t)-\bar X(t)-X_1^\e(t)|^k\\
\ds\les K\sup_{t\in[0,T]}\dbE\Big[\(\int_0^t\big(|\d b(t,s)|+|\D^1_\e b(t,s)|\,|X^\e(s)-\bar X(s)|\big){\bf1}_{[\t,\t+\e]}(s)ds\)^k\\
\ds\qq\qq+\(\int_0^t|\D^1_\e\si(t,s)|^2|X^\e(s)-\bar X(s)|^2{\bf1}_{[\t,\t+\e]}(s)ds\)^{k\over2}\Big]\\
\ns\ds\les K\dbE\Big[\(\int_\t^{\t+\e}\big(1+|\bar X(s)|+|X^\e(s)-\bar X(s)|\big)ds\)^k+\(\int_\t^{\t+\e}|X^\e(s)-\bar X(s)|^2ds\)^{k\over2}\Big]\\
\ds\les K\Big[\e^{k-1}\2n\int_\t^{\t+\e}\3n\big(1+\dbE|\bar X(s)|^k+\dbE|X^\e(s)-\bar X(s)|^k\big)ds+\e^{k-2\over2}\2n\int_\t^{\t+\e}\3n\dbE|X^\e(s)-\bar X(s)|^kds\Big]\\
\ss\ds\les\1n K\Big[\e^k\big(1\1n+\3n\sup_{t\in[0,T]}\2n\dbE|\bar X(t)|^k\2n+\3n\sup_{t\in[0,T]}\dbE|X^\e(t)\1n-\1n\bar X(t)|^k\big)\1n+\1n\e^{k\over2}\2n\sup_{t\in[0,T]}\dbE|X^\e(t)\1n-\1n\bar X(t)|^k\Big]\2n\les\1n K\e^k.\ea\ee
This proves the third relation in \rf{|X_1|,|X_2|}.

\ms

\it Step 3. Estimation of $X^\e(\cd)-\bar X(\cd)-X_1^\e(\cd)-X_2^\e(\cd)$. \rm

\ms

For this step, we first let $X_2^\e(\cd)$ be the solution of \rf{X_2-equation*}. Then for $1\les k\les {p\over2}$,
$$\ba{ll}
\ds\sup_{t\in[0,T]}\1n\dbE|X_2^\e(t)|^k\1n\\
\ds \les \1nK\dbE\Big[\(\1n\int_0^T\3n|X_1^\e(s)|^2ds\)^k\2n+ \1n\(\int_\t^{\t+\e}\3n\3n\big[1+|\bar X(s)|\big]ds\)^k\3n+ \(\1n\int_0^T\3n|X_1^\e(s)|^4
ds\)^{k\over2}\2n+\1n\(\int_\t^{\t+\e}\3n\3n|X_1^\e(s)|^2ds\)^{k\over2}\Big]\\ [3mm]
\ss\ds\les\2n K\dbE\Big[\1n\int_0^T\3n|X_1^\e(s)|^{2k}ds\1n+\1n\e^k\2n+\1n\e^{k-1}\2n
\int_\t^{\t+\e}\3n\3n|X_1^\e(s)|^kds\1n+\2n
\int_0^T\3n|X_1^\e(s)|^{2k}ds\1n+\1n\e^{k-2\over2}\2n\int_\t^{\t+\e}\3n\2n
|X_1^\e(s)|^kds\Big]\1n\les
\2n K\e^k.\ea$$
This proves the forth relation in \rf{|X_1|,|X_2|}. Next, we need to have a further expansion for the coefficients of the state equation. For $1\les i\les n$, we observe
$$\ba{ll}
\ds\wt b_x^i(t,s)\equiv\int_0^1b^i_x\big(t,s,\bar X(s)+\a[X^\e(s)-\bar X(s)],u^\e(s)\big)d\a\\
\ss\ds=\int_0^1\(b^i_x\big(t,s,\bar X(s)+\a[X^\e(s)-\bar X(s)],u^\e(s)\big)-b^i_x(t,s,\bar X(s),u^\e(s))\)d\a\\
\ss\ds\qq+b^i_x(t,s,\bar X(s),\bar u(s))+\big[b^i_x(t,s,\bar X(s),u)-b^i_x(t,s,\bar X(s),\bar u(s))\big]{\bf1}_{[\t,\t+\e]}(s)\\
\ss\ds=\(\int_0^1\int_0^1\a b^i_{xx}\big(t,s,\bar X(s)+\a\b[X^\e(s)-\bar X(s)],u^\e(s)\big)d\b d\a\)[X^\e(s)-\bar X(s)]\\%
\ss\qq+b_x^i(t,s)+\d b^i_x(t,s){\bf1}_{[\t,\t+\e]}(s)\\
\ss\ds=\Big[{1\over2}b^i_{xx}(t,s)+{1\over2}\d b^i_{xx}(t,s){\bf1}_{[\t,\t+\e]}(s)+\D^2_\e b^i_{xx}(t,s)\Big][X^\e(s)-\bar X(s)]\\
\ss\qq+b_x^i(t,s)+\d b^i_x(t,s){\bf1}_{[\t,\t+\e]}(s),\ea$$
where $b_{xx}^i(t,s)$ and $\d b_x(t,s)$ are as in \rf{abb}, and
$$\ba{ll}
\ns\ds \D^2_\e b^i_{xx}(t,s)\1n=\2n\int_0^1\3n\int_0^1\3n\a\Big[b^i_{xx}\big(t,s,\bar X(s)\1n+\1n\a\b[X^\e(s)\1n-\1n\bar X(s)],u^\e(s)\big)\1n-\1n b^i_{xx}(t,s,\bar X(s),u^\e(s))\Big]d\b d\a.
\ea$$
%
%
Combining the above, we see that
\bel{b}\ba{ll}
\ss\ds b(t,s,X^\e(s),u^\e(s))-b(t,s,\bar X(s),\bar u(s))\\
\ss\ds=b_x(t,s)[X^\e(s)-\bar X(s)]+{1\over2}b_{xx}(t,s)[X^\e(s)-\bar X(s)]^2+\d b(t,s){\bf1}_{[\t,\t+\e]}(s)\\
\ss\ds\qq+\(\d b_x(t,s)[X^\e(s)-\bar X(s)]+{1\over2}\d b_{xx}(t,s)
[X^\e(s)-\bar X(s)]^2\){\bf1}_{[\t,\t+\e]}(s)\\
\ss\ds\qq+\D_\e^2b_{xx}(t,s)[X^\e(s)-\bar X(s)]^2\\
\ss\ds\equiv b_x(t,s)[X^\e(s)-\bar X(s)]+{1\over2}b_{xx}(t,s)[X^\e(s)-\bar X(s)]^2+\d b(t,s){\bf1}_{[\t,\t+\e]}(s)+R_\e^b(t,s),\ea\ee
where
$$b_{xx}(t,s)[X^\e(s)-\bar X(s)]^2=\begin{pmatrix}[X^\e(s)-\bar X(s)]^\top b^1_{xx}(t,s)[X^\e(s)-\bar X(s)]\\
\vdots\\ [X^\e(s)-\bar X(s)]^\top b^n_{xx}(t,s)[X^\e(s)-\bar X(s)]\end{pmatrix}.$$
The terms $\d b_{xx}(t,s)[X^\e(s)-\bar X(s)]^2$ and $\D_\e^2b_{xx}(t,s)[X^\e(s)-\bar X(s)]^2$ are defined similarly, and we have denoted
\bel{Rb}\ba{ll}
\ss\ds R_\e^b(t,s)=\(\d b_x(t,s)[X^\e(s)-\bar X(s)]+{1\over2}\d b_{xx}(t,s)
[X^\e(s)-\bar X(s)]^2\){\bf1}_{[\t,\t+\e]}(s)\\
\ss\ds\qq\qq\qq+\D_\e^2b_{xx}(t,s)[X^\e(s)-\bar X(s)]^2.\ea\ee
Likewise, we have
\bel{si}\ba{ll}
\ss\ds\si(t,s,X^\e(s),u^\e(s))-\si(t,s,\bar X(s),\bar u(s))\\
%
%
%
%
\ss\ds=\si_x(t,s)[X^\e(s)-\bar X(s)]+{1\over2}\si_{xx}(t,s)[X^\e(s)-\bar X(s)]^2\\
\ss\ds\qq+\(\d\si(t,s)+\d\si_x(t,s)[X^\e(s)-\bar X(s)]\){\bf1}_{[\t,\t+\e]}(s)+R^\si_\e(t,s),\ea\ee
where
\bel{Rsi}\ba{ll}
\ss\ds R^\si_\e(t,s)={1\over2}\d\si_{xx}(t,s)
[X^\e(s)-\bar X(s)]^2{\bf1}_{[\t,\t+\e]}(s)+\D_\e^2\si_{xx}(t,s)[X^\e(s)-\bar X(s)]^2,\ea\ee
with $\si_{xx}(t,s)[X^\e(s)-\bar X(s)]^2$ and so on being similarly defined as above. Therefore, we may rewrite \rf{X-X-X1} as follows
\bel{X-X2}\ba{ll}
\ds X^\e(t)-\bar X(t)-X_1^\e(t)\\
\ds=\int_0^t\(b_x(t,s)[X^\e(s)-\bar X(s)-X_1^\e(s)]+{1\over2}b_{xx}(t,s)[X^\e(s)-\bar X(s)]^2\\
\ds\qq\qq+\d b(t,s){\bf1}_{[\t,\t+\e]}(s)+R_\e^b(t,s)\)ds\\
\ds\qq+\int_0^t\(\si_x(t,s)[X^\e(s)-\bar X(s)-X_1^\e(s)]+{1\over2}\si_{xx}(t,s)[X^\e(s)-\bar X(s)]^2\\
\ds\qq\qq+\d\si_x(t,s)[X^\e(s)-\bar X(s)]{\bf1}_{[\t,\t+\e]}(s)\1n+\1n R^\si_\e(t,s)\)dW(s).\ea\ee
Therefore,
$$\ba{ll}
\ss\ds X^\e(t)-\bar X(t)-X_1^\e(t)-X_2^\e(t)=\int_0^t\Big[b_x(t,s)[X^\e(s)-\bar X(s)-X_1^\e(s)-X_2^\e(s)]\\
\ss\ds\qq+{1\over2}b_{xx}(t,s)\([X^\e(s)-\bar X(s)]^2-X_1^\e(s)^2\)+R^b_\e(t,s)\Big]ds\\
\ss\ds\qq+\int_0^t\Big[\si_x(t,s)[X^\e(s)-\bar X(s)-X_1^\e(s)-X_2^\e(s)]+{1\over2}\si_{xx}(t,s)\([X^\e(s)-\bar X(s)]^2-X_1^\e(s)^2\)\\
\ss\ds\qq+\d\si_x(t,s)[X^\e(s)-\bar X(s)-X^\e_1(s)]{\bf1}_{[\t,\t+\e]}(s)+R^\si_\e(t,s)\Big]dW(s).\ea$$
Consequently, by the stability estimate of FSVIEs, we have (for $2\les k\les{p\over2}$)
\bel{X-X-X-X1}\ba{ll}
\ds\sup_{t\in[0,T]}\dbE|X(t)-\bar X(t)-X_1^\e(t)-X_2^\e(t)|^k\\
\ds\les
K\sup_{t\in[0,T]}\dbE\Big[\(\int_0^t\Big|{1\over2}b_{xx}(t,s)
\big([X^\e(s)-\bar X(s)]^2-X_1^\e(s)^2\big)+R^b_\e(t,s)\Big|ds\)^k\\
\ds\qq\qq\q+\(\int_0^t\Big|{1\over2}\si_{xx}(t,s)\big([X^\e(s)-\bar X(s)]^2-X_1^\e(s)^2\big)\\
\ds\qq\qq\qq+\d\si_x(t,s)[X^\e(s)-\bar X(s)-X^\e_1(s)]{\bf1}_{[\t,\t+\e]}(s)+R^\si_\e(t,s)\Big|^2ds\)^{k\over2}\Big]\\
\ds\les\1n K\2n\sup_{t\in[0,T]}\1n\dbE\Big\{\Big[\int_0^t\3n\(|X^\e(s)\1n-\1n\bar X(s)-X_1^\e(s)|^k|X^\e(s)-\bar X(s)+X^\e_1(s)|^k+|R^b_\e(t,s)|^k\)ds\Big]\\
\ds\qq\qq\qq\qq+\Big[\int_0^t\(|X(s)-\bar X(s)-X_1^\e(s)|^k|X^\e(s)-\bar X(s)+X^\e_1(s)|^k\\
\ds\qq\qq\qq\qq\qq+|X^\e(s)-\bar X(s)-X^\e_1(s)|^k{\bf1}_{[\t,\t+\e]}(s)+|R^\si_\e(t,s)|^k\)ds\Big]\Big\}\\
\ds\les K\Big\{\sup_{s\in[0,T]}\dbE\Big[|X^\e(s)-\bar X(s)-X_1^\e(s)|^k\(|X^\e(s)-\bar X(s)|^k+|X^\e_1(s)|^k\)\Big]\\
\ns\ds\qq+\e\sup_{s\in[0,T]}\dbE|X^\e(s)-\bar X(s)-X_1^\e(s)|^k\\
\ds\qq+\sup_{t\in]0,T]}\dbE\int_0^t|R^b_\e(t,s)|^kds+\sup_{t\in]0,T]
}\dbE\int_0^t|R^\si_\e(t,s)|^kds\Big\}\\
\ns\ds\les K\Big\{\e^{3k\over2}+\e^{k+1}+\sup_{t\in]0,T]}\dbE\int_0^t|R^b_\e(t,s)|^kds+\sup_{t\in[0,T]
}\dbE\int_0^t|R^\si_\e(t,s)|^kds\Big\}.\ea\ee
For the last two terms on the right-hand side of above, we note that
$$\ba{ll}
\ns\ds\sup_{t\in[0,T]}\dbE\int_0^t|R^b_\e(t,s)|^kds\\
\ns\ds=\sup_{t\in[0,T]}\dbE\int_0^t\Big|\(\d b_x(t,s)[X^\e(s)-\bar X(s)]+{1\over2}\d b_{xx}(t,s)
[X^\e(s)-\bar X(s)]^2\){\bf1}_{[\t,\t+\e]}(s)\\
\ss\ds\qq\qq\qq\qq\qq+\D_\e^2b_{xx}(t,s)[X^\e(s)-\bar X(s)]^2\Big|^kds\\
\ns\ds\les K\Big[\dbE\int_0^T\(|X^\e(s)-\bar X(s)|^k+|X^\e(s)-\bar X(s)|^{2k}\){\bf1}_{[\t,\t+\e]}(s)ds\\
\ss\ds\qq\qq\qq\qq\qq+\sup_{t\in[0,T]}\dbE\int_0^t|\D^2_\e b_{xx}(t,s)|^k|X^\e(s)-\bar X(s)|^{2k}ds\Big]\\
\ns\ds\les K\Big[\e^{k+2\over2}+\e^k\int_0^t\(\dbE\big[K\land\rho\big(|X^\e(s)-\bar X(s)|\big)\big]^{2k}\)^{1\over2}ds\Big]=K\big[\e^{k+2\over2}+o(\e^k)\big].\ea$$
In the above, $\rho(\cd)$ is the modulus of continiuity for the map $x\mapsto b_{xx}(t,s,x,u)$. Similarly,
$$\ba{ll}
\ns\ds\sup_{t\in[0,T]}\dbE\int_0^t|R^\si_\e(t,s)|^kds\\
\ns\ds=\sup_{t\in[0,T]}\dbE\int_0^t\Big|\({1\over2}\d\si_{xx}(t,s)
[X^\e(s)-\bar X(s)]^2\){\bf1}_{[\t,\t+\e]}(s)+\D_\e^2\si_{xx}(t,s)[X^\e(s)-\bar X(s)]^2\Big|^kds\\
\ns\ds\les K\Big[\dbE\int_0^T\(|X^\e(s)-\bar X(s)|^{2k}\){\bf1}_{[\t,\t+\e]}(s)ds+\sup_{t\in[0,T]}\dbE\int_0^t|\D^2_\e\si_{xx}(t,s)|^k|X^\e(s)-\bar X(s)|^{2k}ds\Big]\\
\ns\ds\les K\Big[\e^{k+1}+\e^k\int_0^t\(\dbE\big[K\land\rho\big(|X^\e(s)-\bar X(s)|\big)\big]^{2k}\)^{1\over2}ds\Big]=K\big[\e^{k+1}+o(\e^k)\big]=o(\e^k).\ea$$
Hence,
\bel{X-X-X-X2}\sup_{t\in[0,T]}\dbE|X(t)-\bar X(t)-X_1^\e(t)-X_2^\e(t)|^k\les K\big[o(\e^k)+\e^{k+2\over2}\big].\ee
For $k=2$, the above gives the last relation in \rf{|X_1|,|X_2|}. It also shows that for $k>2$, the above might not gives $o(\e^k)$.

\ms

\it Step 4. Variation of the cost functional. \rm

\ms

We now look at the variation of the cost functional. Note that
\bel{X_1+X_2}\ba{ll}
\ss\ds X_1^\e(t)+X_2^\e(t)\\
\ss\ds=\int_0^tb_x(t,s)[X^\e_1(s)+X^\e_2(s)]ds+\int_0^t\si_x(t,s)[X^\e_1(s)+X^\e_2(s)]dW(s)\\
\ss\ds\q+\int_0^t\({1\over2}b_{xx}(t,s)X_1^\e(s)^2+\d b(t,s){\bf1}_{[\t,\t+\e]}(s)\)ds\\
\ss\ds\q+\2n\int_0^t\3n\({1\over2}\si_{xx}(t,s)X_1^\e(s)^2+\big[\d\si(t,s)+\d\si_x(t,s)X_1^\e(s)\big]
{\bf1}_{[\t,\t+\e]}(s)\)dW(s),
\q t\in[0,T].\ea\ee
Next we let $(h_x^\top,\z(\cd))$ satisfy the first equation in \rf{1st-adjoint-equation}, and $(Y(\cd),Z(\cd\,,\cd))$ be the adapted M-solution of the second equation
in \rf{1st-adjoint-equation}, for which
$$Y(t)=\dbE Y(t)+\int_0^tZ(t,s)dW(s),\qq(t,s)\in\D_*[0,T].$$
Then by the duality between \rf{X_1+X_2} and the second equation in \rf{1st-adjoint-equation}, as in \cite[Theorem 5.1]{Yong-2008}, one has
$$\ba{ll}
\ds\dbE\int_0^T\lan X_1^\e(t)+X_2^\e(t),g_x(t)^\top+b_x(T,t)^\top h_x^\top+\si_x(T,t)^\top\z(t)\ran dt\\
\ds=\dbE\int_0^T\lan Y(t),\int_0^t\({1\over2}b_{xx}(t,s)X_1^\e(s)^2+\d b(t,s){\bf1}_{[\t,\t+\e]}(s)\)ds\\
\ds\q+\int_0^t\({1\over2}\si_{xx}(t,s)X_1^\e(s)^2+\big[\d\si(t,s)+\d\si_x(t,s)X_1^\e(s)\big]
{\bf1}_{[\t,\t+\e]}(s)\)dW(s)\ran dt\\
\ss\ds=\dbE\int_0^T\int_s^T\lan Y(t),\({1\over2}b_{xx}(t,s)X_1^\e(s)^2+\d b(t,s){\bf1}_{[\t,\t+\e]}(s)\)\ran dtds\\
\ss\ds\q+\dbE\int_0^T\int_s^T\lan Z(t,s),\({1\over2}\si_{xx}(t,s)X_1^\e(s)^2+\big[\d\si(t,s)+\d\si_x(t,s)X_1^\e(s)\big]
{\bf1}_{[\t,\t+\e]}(s)\)\ran dtds.\ea$$
Also, by \rf{X_1+X_2} and the first equation in \rf{1st-adjoint-equation},
$$\ba{ll}
\ds\dbE\(h_x[X_1^\e(T)+X_2^\e(T)]\)\\
\ss\ds=\dbE\Big[\int_0^Th_x\(b_x(T,s)[X_1^\e(s)+X_2^\e(s)]+{1\over2}b_{xx}(T,s)X_1^\e(s)^2+\d b(T,s){\bf1}_{[\t,\t+\e]}(s)\)ds\\
\ds\q+\int_0^T\z(s)^\top\(\si_x(T,s)[X_1^\e(s)+X_2^\e(s)]+{1\over2}\si_{xx}(T,s)X_1^\e(s)^2\\
\ds\qq\qq\qq+\big(\d\si(T,s)+\d\si_x(T,s)X_1^\e(s)\big){\bf1}_{[\t,\t+\e]}(s)\)ds\Big]\\
\ds=\dbE\int_0^T\Big[\lan b_x(T,s)^\top h_x^\top+\si_x(T,s)^\top\z(s),X_1^\e(s)+X_2^\e(s)\ran\\
\ss\ds\qq\qq+{1\over2}h_xb_{xx}(T,s)X_1^\e(s)^2+{1\over2}\z(s)^\top\si_{xx}(T,s)X_1^\e(s)^2\\
\ss\ds\qq\qq+\(h_x\d b(T,s)+\z(s)^\top\big(\d\si(T,s)+\d\si_x(T,s)X_1^\e(s)\big)\){\bf1}_{[\t,\t+\e]}(s)\Big]ds.\ea$$
Thus,
$$\ba{ll}
\ds\dbE\(h_x[X_1^\e(T)+X_2^\e(T)]+\int_0^Tg_x(t)[X_1^\e(t)+X_2^\e(t)]dt\)\\
\ss\ds=\dbE\int_0^T\int_s^T\lan Y(t),\({1\over2}b_{xx}(t,s)X_1^\e(s)^2+\d b(t,s){\bf1}_{[\t,\t+\e]}(s)\)\ran dtds\\
\ds\q+\dbE\int_0^T\int_s^T\lan Z(t,s),\({1\over2}\si_{xx}(t,s)X_1^\e(s)^2+\big[\d\si(t,s)+\d\si_x(t,s)X_1^\e(s)\big]
{\bf1}_{[\t,\t+\e]}(s)\)\ran dtds\\
\ss\ds\qq\qq+\int_0^T\Big[{1\over2}h_xb_{xx}(T,s)X_1^\e(s)^2+{1\over2}\z(s)^\top\si_{xx}(T,s)X_1^\e(s)^2\\
\ss\ds\qq\qq+\(h_x\d b(T,s)+\z(s)^\top\big(\d\si(T,s)+\d\si_x(T,s)X_1^\e(s)\big)\){\bf1}_{[\t,\t+\e]}(s)\Big]ds\\
\ss\ds=\dbE\Big\{\int_0^T\int_s^T\Big[{1\over2}\(\lan Y(t),b_{xx}(t,s)X_1^\e(s)^2\ran
+\lan Z(t,s),\si_{xx}(t,s)X_1^\e(s)^2\ran\)\\
\ss\ds\qq\qq+\(\lan Y(t),\d b(t,s)\ran+\lan Z(t,s),
\d\si(t,s)+\d\si_x(t,s)X_1^\e(s)\ran\){\bf1}_{[\t,\t+\e]}(s)dtds\Big]\\
\ss\ds\qq\qq+\int_0^T\Big[{1\over2}\(\lan h_x^\top,b_{xx}(T,s)X_1^\e(s)^2\ran+\lan\z(s),\si_{xx}(T,s)X_1^\e(s)^2\ran\)\\
\ss\ds\qq\qq+\(\lan h_x^\top,\d b(T,s)\ran+\lan\z(s),\d\si(T,s)+\d\si_x(T,s)X_1^\e(s)\ran\) {\bf1}_{[\t,\t+\e]}(s)\Big]ds\Big\}\\
\ss\ds=\dbE\Big\{\int_0^T{1\over2}\Big[\lan h_x^\top,b_{xx}(T,s)X_1^\e(s)^2\ran+\lan\z(s),\si_{xx}(T,s)X_1^\e(s)^2\ran\\
\ds\qq\qq+\int_s^T\(\lan Y(t),b_{xx}(t,s)X_1^\e(s)^2\ran+\lan Z(t,s),\si_{xx}(t,s)X_1^\e(s)^2\ran\)dt\Big]ds\\
\ds\qq\qq+\int_0^T\Big[ \lan h_x^\top,\d b(T,s)\ran+\lan\z(s),
\d\si(T,s)\ran\\
\ds\qq\qq+\int_s^T\(\lan Y(t),\d b(t,s)\ran+\lan Z(t,s),
\d\si(t,s)\ran\)dt\Big]{\bf1}_{[\t,\t+\e]}(s)ds\\
\ss\ds\qq\qq+\int_0^T\(\lan\z(s),\d\si_x(T,s)X_1^\e(s)\ran+\int_s^T\lan Z(t,s),\d\si_x(t,s)X_1^\e(s)\ran dt\){\bf1}_{[\t,\t+\e]}(s)ds\Big\}.\ea$$
Note that
$$\ba{ll}
\ds\dbE\Big|\int_0^T\(\lan\z(s),\d\si_x(T,s)X_1^\e(s)\ran+\int_s^T\lan Z(t,s),\d\si_x(t,s)X_1^\e(s)\ran dt\){\bf1}_{[\t,\t+\e]}(s)ds\Big|\\
\ds\les K\int_\t^{\t+\e}\dbE\Big[\(|\z(s)|+\int_s^T|Z(t,s)|dt\)|X_1^\e(s)|\Big]ds\\
%
%
\ds\les K\Big[\int_\t^{\t+\e}\dbE\(|\z(s)|+\int_s^T|Z(t,s)|dt\)^2ds\Big]^{1\over2}
\Big[\int_\t^{\t+\e}\dbE|X_1^\e(s)|^2ds\Big]^{1\over2}\\
\ds\les K\e\Big[\int_\t^{\t+\e}\dbE\(|\z(s)|+\int_s^T|Z(t,s)|dt\)^2ds\Big]^{1\over2}=o(\e).\ea$$
Now, similar to the above Steps 1--3, we have
%
$$\ba{ll}
\ds J(u^\e(\cd)-J(\bar u(\cd))\\
\ds=\dbE\Big[h(X^\e(T))-h(\bar X(T))+\int_0^T\(g(s,X^\e(s),u^\e(s))-g(s,\bar X(s),\bar u(s))\)ds\Big]\\
\ss\ds=\dbE\Big[h_x[X^\e(T)-\bar X(T)]+{1\over2}h_{xx}[X^\e(T)-\bar X(T)]^2+\D_\e^2h_{xx}(T)[X^\e(T)-\bar X(T)]^2\\
\ss\ds\qq+\int_0^T\(g_x(s)[X^\e(s)-\bar X(s)]+{1\over2}g_{xx}(s)[X^\e(s)-\bar X(s)]^2\\
\ss\ds\qq+\big(\d g(s)+\d g_x(s)[X^\e(s)-\bar X(s)]+{1\over2}\d g_{xx}(s)
[X^\e(s)-\bar X(s)]^2\big)I_{[\t,\t+\e]}(s)\\
\ss\ds\qq+\D_\e^2g_{xx}(s)[X^\e(s)-\bar X(s)]^2\)ds\Big] \\
\ss\ds=\dbE\Big[h_x[X_1^e(T)+X_2^\e(T)]+\2n\int_0^T\3n g_x(s)[X_1^\e(s)\1n+\1nX_2^\e(s)]ds\\
\ss\ds\q+{1\over2}X_1^\e(T)^\top h_{xx}X_1^\e(T)+\int_0^T\({1\over2}X_1^\e(s)^\top g_{xx}(s)X_1^\e(s)\1n+\1n\d g(s){\bf1}_{[\t,\t+\e]}(s)\)ds\Big]\1n+\1n o(\e)\\
\ss\ds={1\over2}\dbE\Big\{h_{xx}(\bar X(T))X_1^\e(T)^2+\int_0^T\Big[\lan h_x(\bar X(T))^\top,b_{xx}(T,s,\bar X(s),\bar u(s))X_1^\e(s)^2\ran\\
\ns\ds\qq+\lan\z(s),\si_{xx}(T,s,\bar X(s),\bar u(s))X_1^\e(s)^2\ran+g_{xx}(s,\bar X(s),\bar u(s))X_1^\e(s)^2\\
\ns\ds\qq+\int_s^T\(\lan Y(t),b_{xx}(t,s,\bar X(s),\bar u(s))X_1^\e(s)^2\ran\\
\ns\ds\qq+\lan Z(t,s),\si_{xx}(t,s,\bar X(s),\bar u(s))X_1^\e(s)^2\ran\)dt\Big]ds\Big\}\ea$$
$$\ba{ll}
\ss\ds\qq\qq+\dbE\int_\t^{\t+\e}\Big[\lan h_x(\bar X(T))^\top,b(T,s,\bar X(s),u)-b(T,s,\bar X(s),\bar u(s))\ran\\
\ns\ds\qq\qq\qq\qq+\lan\z(s),\si(T,s,\bar X(s),u)-\si(T,s,\bar X(s),\bar u(s))\ran\\
\ns\ds\qq\qq\qq\qq+g(s,\bar X(s),u)-g(s,\bar X(s),\bar u(s))\\
\ss\ds\qq\qq\qq+\int_s^T\(\lan Y(t),b(t,s,\bar X(s),u)-b(t,s,\bar X(s),\bar u(s))\ran\\
\ns\ds\qq\qq\qq+\lan Z(t,s),
\si(t,s,\bar X(s),u)-\si(t,s,\bar X(s),\bar u(s))\ran\)dt\Big]ds+o(\e).\ea$$
By defining the Hamiltonian $H(\cd)$ as in \rf{H}, we obtain \rf{J-J}. \endpf

\section{Concluding Remarks}

In this paper, we have developed a spike variation technique for optimal controls of FSVIEs. The main contribution is the derivation of the second-order adjoint equation which is a little different from a standard BSVIE. Thus the well-posedness of such a BSVIE is a part of novelty. We have used a trick used in the derivation of Riccatti equation in (stochastic) linear-quadratic optimal control, together with the introduction of the auxiliary process so that the It\^o's formula can be used. 

\ms

As an extension, we have seen that for multi-person dynamic games of FSVIEs, one can obtain the corresponding maximum principle for the Nash equilibria. It is natural to ask what happens when $N\to\i$? Under certain structure conditions, this will be related to the mean-field dynamic games of FSVIEs. The relevant results will be reported in our forthcoming work in the near future.

\end{document}